\documentclass[sigconf,natbib=True,anonymous=False,authorversion,nonacm]{acmart}

\usepackage{algorithm}
\usepackage{algorithmic}

\usepackage[utf8]{inputenc} 
\usepackage[T1]{fontenc}    
\usepackage{hyperref}       
\usepackage{url}            
\usepackage{booktabs}       
\usepackage{amsfonts}       
\usepackage{nicefrac}       
\usepackage{microtype}      
\usepackage{xcolor}         
\usepackage{enumitem}
\usepackage{dsfont}
\usepackage{hyperref}

\newtheorem*{theorem*}{Theorem}

\newtheorem{theorem}{Theorem}[section]
\newtheorem{lemma}{Lemma}[section]
\newtheorem{observation}{Observation}[section]

\newtheorem{assumption}{Assumption}[section]

\DeclareMathOperator*{\argmin}{arg\,min}
\DeclareMathOperator*{\tr}{\text{Tr}}

\setcopyright{none}

\begin{document}

\title{Momentum-inspired Low-Rank Coordinate Descent for Diagonally Constrained SDPs}









\author{Junhyung Lyle Kim, Jose Antonio Lara Benitez, Mohammad Taha Toghani, Cameron Wolfe, Zhiwei Zhang, Anastasios Kyrillidis}
\thanks{Authors before AK are listed in alphabetical order.}
\affiliation{%
  \institution{Rice University}
  \streetaddress{6100 Main St.}
  \city{Houston,~TX}
  \country{USA}}
\email{jlylekim, jx24, mt72, crw13, zhiwei, anastasios@rice.edu}

\renewcommand{\shortauthors}{Kim, et al.}

\begin{abstract}
We present a novel, practical, and provable approach for solving diagonally constrained semi-definite programming (SDP) problems at scale \emph{using accelerated non-convex programming}.  
Our algorithm non-trivially combines acceleration motions from convex optimization with coordinate power iteration and matrix factorization techniques.
The algorithm is extremely simple to implement, and adds only a single extra hyperparameter -- momentum. 
We prove that our method admits local linear convergence in the neighborhood of the optimum and always converges to a first-order critical point.
Experimentally, we showcase the merits of our method on three major application domains: MaxCut, MaxSAT, and MIMO signal detection. 
In all cases, our methodology provides significant speedups over non-convex and convex SDP solvers -- $5\times$ faster than state-of-the-art non-convex solvers, and $9$ to $10^3 \times$ faster than convex SDP solvers -- with comparable or improved solution quality.
\end{abstract}



\maketitle

\section{Introduction}\label{sec:intro}

\textbf{Background.}
This work focuses on efficient ways to solve semidefinite programming instances (SDPs) with diagonal constraints:
\begin{equation}
\begin{aligned}
\min_{X \succeq 0} \quad & \left \langle C, X \right \rangle  & \textrm{s.t.} \quad & X_{ii} = 1, ~i = 1, \dots, n.
\end{aligned} \label{eq:continuous_form}
\end{equation}
Here, $X \in \mathbb{S}^n$ is the $n \times n$ symmetric optimization variable and $C\in \mathbb{S}^n$ is a problem-dependent cost matrix.
The above formulation usually appears in practice as the convex relaxation of quadratic form optimization over discrete variables.
\begin{equation}
\begin{aligned}
& \underset{x \in \mathcal{D}^n}{\min}
& & \langle x, C x \rangle,
\end{aligned} \label{eq:discrete_form}
\end{equation}
$\mathcal{D}$ is a discrete set on the unit imaginary circle:
\begin{align*}
\mathcal{D} = \left\{ e^{\frac{\imath 2\pi m}{M}} ~:~ m = 0, 1, \dots, M-1 \right\}.
\end{align*}
where $\imath := \sqrt{-1}$ denotes the imaginary unit.
Different values for $M$ define different realizations of the discrete set $\mathcal{D}$. 
For example, $M=2$ reduces $\mathcal{D}$ to the binary set $\{ \pm 1 \}$, while $\mathcal{D}$ becomes $\{\pm 1, \pm \imath\}$ for $M = 4$.\footnote{Many other realizations of the set $\mathcal{D}$ exist. 
For example, we define \mbox{$\mathcal{D} = \{1, e^{\frac{\imath 2\pi}{3}}, e^{\frac{\imath 4\pi}{3}} \}$} for Max-$3$-Cut~\cite{goemans2004approximation,so2007approximating}.}
The formulation in~\eqref{eq:discrete_form} appears in many applications, including
MaxCut or Max-$k$-Cut~\cite{karp1972reducibility,barahona1988application,goemans1995improved,deza1994applications,hartmann1996cluster,shi2000normalized,mei2017solving,frieze1995improved, goemans2004approximation, so2007approximating}, 
MaxSAT~\cite{goemans1995improved},
maximal margin classification~\cite{gieseke2013polynomial}, 
semi-supervised learning~\cite{wang2013semi},
correlation clustering on dot-product graphs~\cite{veldt2017correlation}, 
community detection~\cite{hajek2016achieving,abbe2018community}, 
and (quantized) phase synchronization in communication networks~\cite{boumal2016nonconvex,zhong2018near}. 
Other problems that can be cast as (quantized) quadratic form optimization over the unit complex torus include phase recovery~\cite{waldspurger2015phase}, angular synchronization~\cite{singer2011angular},
and optimization problems in communication systems~\cite{heath1998simple,love2003grassmannian,motedayen2006optimal,kyrillidis2014fixed,kyrillidis2011rank}.

Typically, solving~\eqref{eq:discrete_form} is computationally expensive due to the \emph{presence of discrete structures}.
Discrete algorithmic solutions have been developed to solve~\eqref{eq:discrete_form} with near-optimal performance and rigorous theoretical guarantees by exploiting the problem's structure~\cite{marti2009advanced, optimization2014inc,mosek2015mosek,krislock2017biqcrunch}.
However, the go-to techniques for solving~\eqref{eq:discrete_form} are oftentimes continuous.
For example, methods exist that relax combinatorial entities in \eqref{eq:discrete_form} and find a solution with non-linear, continuous optimization. 
Such techniques can be roughly categorized into convex \cite{grotschel2012geometric,gartner2012approximation,nesterov1989self,nesterov1994interior,tran2014inexact,tran2016single} and non-convex, continuous  approaches~\cite{burer2003nonlinear,bhojanapalli2018smoothed,bhojanapalli2016dropping,kyrillidis2018provable,wang2017mixing}.

Recently, numerous algorithms have been proposed for solving large-scale SDP instances of the form \eqref{eq:continuous_form}.
\cite{wang2017mixing} present a low-rank coordinate descent approach to solving MaxCut and MaxSat SDP instances that is 10-100$\times$ faster than state-of-the-art solvers.
Similarly,~\cite{yurtsever2019scalable} solve MaxCut instances on a laptop with almost $8$ million vertices and $10^{13}$ matrix entries, while~\cite{goto2019combinatorial} propose a simulated bifurcation algorithm that can solve all-to-all, $100,000$-node MaxCut problems with continuous weights.
Despite such developments, however, SDPs are still hard problems to solve in practice.
\emph{This work aims to contribute along the path of efficient algorithms for solving large-scale SDPs.}

\vspace{0.2cm}
\noindent \textbf{This Paper.}
We present an algorithmic prototype for solving diagonally constrained SDPs at scale.
Based on the low $k$-rank property of SDPs~\cite{alizadeh1997complementarity,barvinok1995problems,pataki1998rank} (see Lemma \ref{lemma:SDPlowrank}),  
we focus on solving a non-convex, but equivalent, formulation of~\eqref{eq:continuous_form}:
\begin{equation}
\begin{aligned}
\min_{V \in \mathbb{R}^{k \times n}}  f(V)\triangleq & \left \langle C, V^\top V \right \rangle  & \textrm{s.t.} \quad & \|v_i\|_2 = 1, ~\forall i.
\end{aligned} \label{eq:noncvx_continuous_form}
\end{equation}
Here, $V \in \mathbb{R}^{k \times n}$ is the optimization variable and the linear constraints on $X$ in~\eqref{eq:continuous_form} translate to $V \in (S^{k-1})^{n}$ (i.e., each column vector $v_i$ satisfies $\|v_i\|_2 = 1$).
Our algorithm, named \texttt{Mixing Method++}, adds momentum to a coordinate power iteration-style technique for solving SDPs by incorporating ideas from accelerated convex optimization and matrix factorization.
As its predecessor~\cite{wang2017mixing}, the algorithm is simple to implement, requiring only one additional hyperparameter -- momentum.
From an empirical perspective, we test the algorithm on MaxCut, MaxSAT, and MIMO signal detection applications, where we demonstrate significant speedups in comparison to other state-of-the-art solvers (e.g., $9\times$ to $10^2\times$ speedup on MaxSAT instances) with comparable or improved solution quality.
Additionally, we prove that our accelerated method admits local linear convergence in the neighborhood of the optimum and always converges to a first order stationary point.
To the best of our knowledge, this is the first theoretical result that incorporates momentum into a coordinate-descent approach for solving SDPs.

\section{Background}

\noindent
\textbf{Notation.} 
We use $\mathbb{S}^n$ to denote the cone of positive semi-definite (PSD) matrices. 
We denote the $k$-dimensional sphere by $S^{k-1}$, where $k-1$ indicates the dimension of the manifold.
The product of $n$, $k$-dimensional spheres is denoted by $(S^{k-1})^{n}$.
For any $v\in \mathbb{R}^n$ and $M\in \mathbb{R}^{n\times m}$, we define $\|v\|_2 $ and $\|M\|_F $ as the $\ell_2$ and Frobenius norms.
For a matrix $V \in \mathbb{R}^{k \times n}$, we denote its $i$-th column interchangeably as $V_{:, i}$ or $v_i$. 
For $a \in \mathbb{R}^n$, we define $D_a$ as the $n\times n$ diagonal matrix with entries $(D_a)_{ii} = a_i$. 
We denote the vector of all $1$'s as $\mathds{1}$ and the Hadamard product as $\odot$.
Singular values are denoted as $\sigma(\cdot)$, while $\sigma_{\texttt{nnz}}(\cdot)$ represents minimum non-zero singular value(s).
The matrix $C$ corresponds to the cost matrix of the optimization problem~\eqref{eq:continuous_form}. 

\medskip
\noindent \textbf{Low-rank Property of SDPs.}
Focusing on~\eqref{eq:continuous_form}, the number of linear constraints on $X$ is far less than the number of variables within $X$. 
This constitutes such SDPs as \emph{weakly constrained}, which yields the following result \cite{barvinok1995problems, pataki1998rank}.

\begin{lemma} \label{lemma:SDPlowrank}
The SDP in~\eqref{eq:continuous_form} has a solution with rank $k = \left \lceil \sqrt{2n}\, \right \rceil$.
\end{lemma}
One can enforce low-rank solutions to weakly constrained SDPs by defining $X = V^\top V$, where $V \in \mathbb{R}^{k \times n}$ and $k(k+1)/2>n$ ~\cite{burer2003nonlinear, boumal2016nonconvex}.
This low-rank parameterization, shown in \eqref{eq:noncvx_continuous_form}, yields a non-convex formulation of \eqref{eq:continuous_form} with the same global solution. 
This formulation can be solved efficiently because the conic constraint in~\eqref{eq:continuous_form} is automatically satisfied.\footnote{Handling a PSD constraint requires an eigenvalue decomposition of $n \times n$ matrices, leading to a $O(n^3)$ overhead per iteration.}

\medskip 
\noindent
\textbf{Mixing Method.}
Consider again the low-rank SDP parameterization given in \eqref{eq:noncvx_continuous_form}.
To connect this formulation with its discrete form in~\eqref{eq:discrete_form}, assume $k = 1$, which yields the following expression:
\begin{equation}
\begin{aligned}
\min_{v \in \mathbb{R}^{1 \times n}} \quad & v C v^\top & \textrm{s.t.} \quad & |v_i| =  1, ~\forall i,
\end{aligned} \label{eq:continuous_form_noncvx2}
\end{equation}
Notice that $v$ is normalized entry-wise during each optimization step.
To solve~\eqref{eq:continuous_form_noncvx2}, one can consider using related algorithms, such as a power iteration-style methods~\cite{mises1929praktische}, as shown below. 
\begin{align*}
    u &\leftarrow Cv, \quad 
    v^{+} \leftarrow \texttt{normalize}(u),
\end{align*}
Here, $\texttt{normalize}(u)$ projects $u$ such that $|v^{+}_i| = 1$ for each entry of $v^{+}$.
This is the crux of the Mixing Method~\cite{wang2017mixing}, which applies a coordinate power iteration routine to sequentially update each column of $V$ and solve~\eqref{eq:noncvx_continuous_form}; see Algorithm~\ref{algo:mixing}.
Variants of Algorithm~\ref{algo:mixing} that utilize random coordinate selection instead of cyclic updates for the columns of $V$ have also been explored \cite{erdogdu2018convergence}.

In~\eqref{eq:noncvx_continuous_form}, diagonal elements of $C$ do not contribute to the optimization due to the constraint $\|v_i\|_2 = 1$, and thus can be set to zero.
When considering only the $i$-th column of $V$, the objective $\left \langle C, V^\top V \right \rangle$ takes the form $2v_i^\top \sum_{j = 1}^n c_{ij} v_j$.
The full gradient of the objective with respect to $v_i$ can be computed in closed form at each iteration, yielding the coordinate power iteration routine outlined in Algorithm~\ref{algo:mixing}.
Interestingly, this procedure also comes with theoretical guarantees~\cite{wang2017mixing}, which we include below for completeness.\footnote{Most of global guarantees in~\cite{wang2017mixing} hold for a variant of the Mixing Method that resembles coordinate gradient descent with step size, rather than the coordinate power iteration method that is used and tested in practice.}
\begin{theorem}
The Mixing Method converges linearly to the global optimum of~\eqref{eq:noncvx_continuous_form}, under the assumption that the initial point is close enough to the optimal solution.
\end{theorem}

\begin{algorithm}[t]
    \caption{Mixing Method~\cite{wang2017mixing}}
    \label{algo:mixing}
\begin{algorithmic}[1]
    \STATE {\bfseries Input:} $C\in \mathbb{S}^n$, $V \in (S^{k-1})^{n}$.
    \WHILE{not yet converged}
        \FOR{$i=1$ {\bfseries to} $n$}
            \STATE $v_i \gets \texttt{normalize}\left(-\sum_{j\neq i} c_{ij}v_j\right)$\label{step:mixing inner update}
        \ENDFOR
    \ENDWHILE
\end{algorithmic}
\end{algorithm}

\section{\texttt{Mixing Method++}}\label{sec:accelerated_mixing}

\noindent
\textbf{Towards a Momentum-inspired Update Rule.}
In this work, we introduce \texttt{Mixing Method++}, which uses acceleration techniques~\cite{polyak1987introduction,nesterov2013introductory} to significantly improve upon the performance of Mixing Method.
For gradient descent, a classical acceleration technique is the Heavy-Ball method~\cite{polyak1987introduction}, which iterates as follows:
\begin{align*}
    w^{t+1} = w^t - \eta \nabla f(w^t) + \beta (w^t - w^{t-1}).
\end{align*}
Here, $w$ is the optimization variable, $f(\cdot)$ is a differentiable loss function, $\eta$ is the step size, and $\beta$ is the momentum parameter. 
Intuitively, the heavy-ball method exploits the history of previous updates by moving towards the direction $w^t - w^{t-1}$, weighted by the momentum parameter $\beta$.
A similar momentum term can be naively incorporated into power iteration: 
\begin{align*}
    w^{t+1} = \texttt{normalize}\left(Cw^{t} + \beta (Cw^{t} - w^{t})\right),
\end{align*}
where there is no notion of a step size $\eta$.
By adapting this accelerated power iteration scheme to the Mixing Method update rule, we arrive at the following recursion:
\begin{align*}
    \hat{v}_i = \texttt{normalize}(g_i + \beta(g_i - v_{i})),
\end{align*}
where $g_i = -\sum_{j \neq i} c_{ij} v_j$.
Because $g_i$ is not normalized before the addition of the momentum term, $g_i$ and $v_i$ can be of significantly different magnitude.
As such, we normalize $g_i$ as an intermediate step and adopt a double-projection approach per iteration; see Algorithm~\ref{algo:accelerated_mixing}.
Such double projection ensures that $u_i$ and $v_i$ are of comparable magnitude when the momentum term is added, resulting in significantly-improved momentum inertia.

\begin{algorithm}[t]
    \caption{\texttt{Mixing Method++}}
    \label{algo:accelerated_mixing}
\begin{algorithmic}[1]
    \STATE {\bfseries Input:} $C\in \mathbb{S}^n$, $V \in (S^{k-1})^{n}$, $0\leq\beta<1$.
    \WHILE{not yet converged}
    \FOR{$i=1$ {\bfseries to} $n$}
        \STATE $u_i \gets \texttt{normalize}\left(-\sum_{j\neq i} c_{ij}v_j\right)$\label{step:mixing++ inner update}
        \STATE $v_i\gets \texttt{normalize}\left(u_i + \beta \left(u_i - v_i\right)\right)$\label{step:mixing++ inner momentum}
    \ENDFOR
    \ENDWHILE
\end{algorithmic}
\end{algorithm}

\medskip \noindent
\textbf{Properties of \texttt{Mixing Method++}.}
The proposed algorithm is comprised of a two-step update procedure.
The last normalization projects \mbox{$u_i+\beta \left(u_i-v_i\right)$} onto the unit sphere, which ensures the updated matrix is feasible. 
In contrast, the first normalization ensures that \mbox{$-\sum_{j=1}^n c_{ij}v_j$} and $v_i$ are comparable in magnitude, as the former may have a norm greater than $1$.
The overall computational cost of the added momentum step is negligible, making the step-wise computational complexity of \texttt{Mixing Method++} roughly equal to that of the Mixing Method.

We show that \texttt{Mixing Method++} converges linearly when $V$ lies in a neighborhood of the global optimum and always reaches a first-order stationary point.
Despite the similarity of our convergence rates to previous work~\cite{wang2017mixing}, achieving acceleration in theory is not always feasible, even in convex, non-stochastic cases; see~\cite{assran2020convergence, devolder2014first,ghadimi2015global,polyak1987introduction,lessard2016analysis,loizou2017momentum}.
Nonetheless, we observe acceleration in practice, and \texttt{Mixing Method++} with $\beta = 0.8$ yields a robust algorithm with significant speed-ups in comparison to Mixing Method.

\section{Convergence Analysis}
\label{sec:convergence}

We begin by presenting some notation for brevity and readability.
$f^\star$ is used to denote the optimal value of the objective function $f(V) = \langle C, V^\top V \rangle$.
The matrices $V$ and $\hat{V}$ refer to the current and next iterates from the inner iteration of Algorithm~\ref{algo:accelerated_mixing}.
For each $i\in [n]$, we define the family of maps $Z_i : \mathbb{R}^{k\times n} \to \mathbb{R} ^{k\times n}$ as: 
\begin{equation}\label{eq:Z_i}
    V\longmapsto \begin{pmatrix}
    \vert  & & \vert & \vert & \vert && \vert \\
    \hat{v_1}   &\ldots & \hat{v}_{i-1} & v_i & v_{i+1} & \ldots & v_n\\
    \vert  & & \vert & \vert & \vert & & \vert
\end{pmatrix}.
\end{equation}
In words, $Z_i(V)$ is the partially updated $V$ matrix (i.e., within the inner loop of Algorithm \ref{algo:accelerated_mixing}) before updating $v_i$ to $\hat{v}_i$.
Thus, only the vector $v_i$ changes from $Z_i(V)$ to $Z_{i+1}(V)$.
Notice that $Z_{n+1}(V) = \hat{V}$. 

We now present our main theoretical results.
All proofs are provided in Section \ref{sec:proofs}.\footnote{The proof that \texttt{Mixing Method++} always converges to a first-order critical point is provided in Section \ref{sec:critical}.}
We adopt Assumption~3.1 from \cite{wang2017mixing}.

\begin{assumption}\label{assu:delta} For $i\in [n]$, assume $\|\sum_{j=1}^n c_{ij}v_j\|_2$ do not degenerate in the procedure. 
That is, all norms are always greater than or equal to a constant value $\delta>0$. 
\end{assumption}

\begin{figure}
    \centering
    \includegraphics[width=0.7\linewidth]{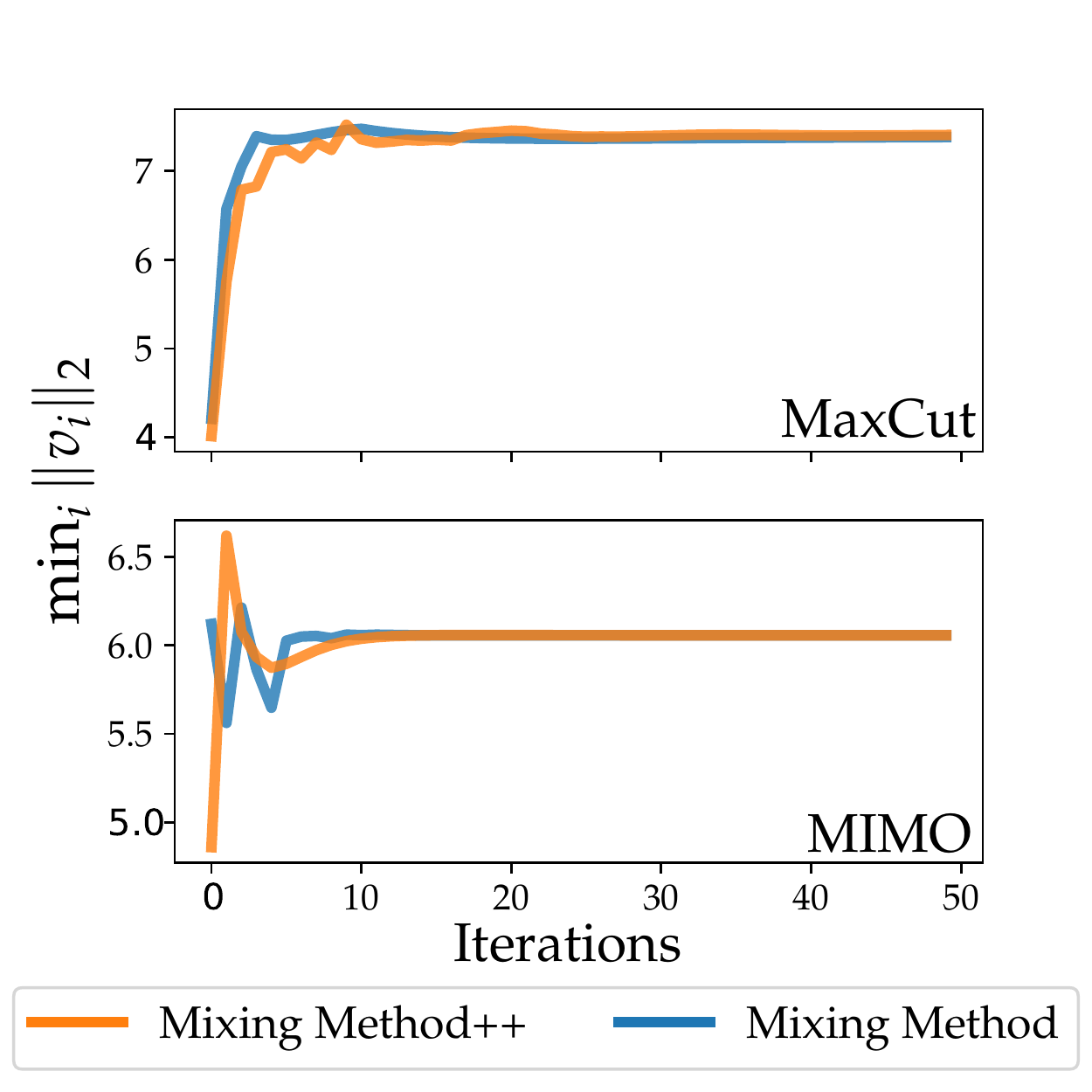}
    \caption{The minimum norm of the columns of $V$ throughout iterations of Mixing Method and \texttt{Mixing Method++} on the MaxCut and MIMO objectives.}
    \label{fig:norm_analysis}
\end{figure}

Put simply, Assumption \ref{assu:delta} states that the update rule will never yield $u_i$ with zero norm.
We find that this assumption consistently holds in practice for both Mixing Method and \texttt{Mixing Method++} applied over several objectives.
In particular, we track the smallest value of $\|\sum_{j \neq i} c_{ij} v_j\|_2$ observed during the inner loop of Mixing Method and Mixing Method++, then plot these minimum norms throughout all iterations in Figure \ref{fig:norm_analysis}.
As can be seen, Assumption \ref{assu:delta} consistently holds for both MaxCut and MIMO objectives.
Furthermore, if Assumption~\ref{assu:delta} does not hold, we can exploit the simple trick from \cite{erdogdu2018convergence} in practice.
Given this assumption, no further assumptions are required for the second normalization step provided that $0\leq \beta < 1$, as summarized below.

\begin{observation}\label{observation0}
  $ w_i = \|u_i-\beta v_i\|_2 \not = 0$, whenever $0\leq\beta< 1$. 
\end{observation}

\begin{proof} Based on upper and lower triangle inequality bounds and the fact that $u_i$ and $v_i$ are both in $S^{k-1}$, we get the inequality \mbox{$1\le 
\left\lVert\left(1+\beta\right)u_i-\beta v_i\right\rVert_2\le 1+2\beta$}. 
\end{proof}

Our first technical lemma lower bounds the decrease in the objective over inner iterations of Algorithm \ref{algo:accelerated_mixing}.

\begin{lemma}
\label{lemma:decreasing} Let $\hat{V}, V$ be the next and previous iteration of Algorithm~\eqref{algo:accelerated_mixing}, respectively.
We define $y_i=\|g_i\|_2$, where $g_i = -\sum_{j \neq i} c_{ij} v_j$.
We then derive the following.
\begin{align*}
     f(V)-f(\hat{V})\ge \frac{1-\beta}{1+\beta}\sum_{i=1}^n y_i\, \|v_i-\hat{v}_i\|_2^2.
\end{align*}
Furthermore, $f$ is non-increasing whenever $0\leq\beta<1$.
\end{lemma}

Next, we lower bound the magnitude of the update from $V$ to $\hat{V}$ based on the objective residual, showing that the magnitude of our update will be large until we approach an optimum.

\begin{lemma}
\label{lemma: lemmad2}
Under assumption~\ref{assu:delta}, there exists positive constants $\gamma$ and $\zeta$ such that
\begin{equation*}
    \| V - \hat{V} \|_F^ 2 \ge \big(\zeta - \gamma\, \|\Delta\|_F  \big)\, (f(V) - f^*),
\end{equation*}
where $\Delta := (D_wD_y+\beta D_y)- (D_{w^*}D_{y^*}+\beta D_{y^*})$.
\end{lemma}

Finally, we upper bound the value of $\|\Delta\|^2_2$ with respect to the objective residual, which later enables the definition of a neighborhood for local linear convergence within Theorem \ref{theorem: linearconvergence}.
\begin{lemma}
\label{lemma: lemmad4}
Under Assumption~\ref{assu:delta}, $\exists~\tau > 0$ such that
\begin{equation*}
    \|(w+\beta \mathds{1})\odot y- (w^*+\beta \mathds{1})\odot  y^*\|_2^2 \le \tau (f(V) - f^*).
\end{equation*}
\end{lemma}

Drawing upon these technical lemmas, we show \texttt{Mixing Method++} converges linearly within a neighborhood of the optimum.

\begin{theorem}[Local linear convergence] \label{theorem: linearconvergence} 
Let $f^\star$ represent the optimal value of the objective function, and let $\delta$ be the non-degenerative lower-bound of Assumption~\ref{assu:delta}. Define the neighborhood,
$$\mathcal{N}_\star : = \left\{ V \in (S^{k-1})^n\, : \,  \tau (f(V)- f^\star) \le \left(\tfrac{\zeta - \kappa}{\gamma} \right)^2  \right\},$$
for some positive constant $\kappa$ and $\zeta$, $\tau$, and $\gamma$ as defined in Lemmas ~\ref{lemma: lemmad2} and~\ref{lemma: lemmad4}.
If $~V\in\mathcal{N}_\star$, we then have:
\begin{equation*}
     f(\hat{V})-f^\star \le \left(1- \rho\right)\,\left(f(V)-f^\star \right),
\end{equation*}
where $\rho \in (0, 1]$.
This result shows that Algorithm~\eqref{algo:accelerated_mixing} converges linearly when $V\in\mathcal{N}_\star$.
\end{theorem}

\section{Experimental Results} \label{sec:exp}
We test \texttt{Mixing Method++} on three well-known applications of~\eqref{eq:continuous_form}: MaxCut, MaxSAT, and MIMO signal detection.
We use the same formulation in \cite{wang2017mixing} to solve MaxCut and MaxSAT, while for MIMO we follow the experimental setup in \cite{liu2017discrete}. 
We implement \texttt{Mixing Method++} in \texttt{C} with $\beta=0.8$.\footnote{We extensively tested our theoretical bound on $\beta$ in experiments. After testing all values in the set $[0,1)$ with an interval of $0.01$, we observe that $\beta = 0.8$ indeed yields the best performance on over 90\% of instances.} 
Each experiment is run on a single core in a homogeneous Linux cluster with 2.63 GHz CPU and 16 GB of RAM.
We compare \texttt{Mixing Method++} with wide variety of solvers depending on the application. For MaxCut, we compare with CGAL \cite{CGAL}, Sedumi \cite{sturm1999using}, SDPNAL+ \cite{yang2015sdpnal}, SDPLR \cite{burer2003nonlinear}, MoSeK \cite{aps2019mosek}, and Mixing Method \cite{wang2017mixing}. For MaxSAT, we compare with Mixing Method and Loandra \cite{loandra}. Lastly, for MIMO signal detection, we compare with Mixing Method.

We often report performance in terms of the objective residual, which is defined as the difference between objective values given by each solver and the lowest objective value achieved by any solver. 
More formally, let $\texttt{Obj}_{ij}$ be the value of the objective function given by optimizer $i$ on SDP instance $j$.
Then, the objective residual of optimizer $i$ on instance $j$ is defined as $|\texttt{Obj}_{ij}-\mathop{\min}\limits_{i'\in O}\texttt{Obj}_{i'j}|$, where $O$ is the set of all solvers used on problem instance $j$.


\subsection{MaxCut}
\begin{table}
\centering
\begin{tabular}{ccccccc}
\toprule
Solver  & & \# SI & &  MLOR & & Acc. ($\times$) \\ 
\cmidrule(lr){1-1} \cmidrule(lr){3-3} \cmidrule(lr){5-5} \cmidrule(lr){7-7} 
Sedumi  & & 60  & & -4.30 & & 360.20 \\ 
MoSeK   & & 74  & & -3.24  & & 348.37  \\ 
SDPNAL+ & & 52  & & \textbf{-4.33} & & 1316.35 \\ 
CGAL &  & 92  &  & 0.93 & & 49.19 \\ 
SDPLR & & 95 &  & -2.89 &  & 9.14  \\ 
Mixing Method  & & 106 & & -2.80 & & 5.26 \\ 
\texttt{Mixing Method++} & & \textbf{111} & & -3.87 & & \textbf{1} \\ 
\bottomrule
\end{tabular}
\caption{\small \textbf{Results on 203 MaxCut instances.} SI stands for solved instances; MLOR stands for median(log(objective residual)). \texttt{Mixing Method++} solved most instances efficiently with good precision.} \label{table:maxcut}
\end{table}

The SDP relaxation of the MaxCut problem is given by
\begin{equation*}
     \underset{\|v_i\|=1,}{\mathrm{max}}\,\frac{1}{2}\sum_{ij}c_{ij}\left(\tfrac{1-v_i^\top v_j}{2}\right),
\end{equation*}
where $C$ is the adjacency matrix.
Thus, we aim to minimize the following objective function:
$$
 f_{\text{MaxCut}} = \left\langle C, \tfrac{\mathds{1}_{n \times n}- V^\top V}{4} \right \rangle ~~\text{s.t.}~~ \|v_i\|_2=1~~\forall i,
$$
where $V \in \mathbb{R}^{k \times n}$ and $\mathds{1}_{n \times n}$ is the $n\times n$ matrix of all $1$s.
We consider the following benchmarks within our MaxCut experiments.

\begin{itemize}[leftmargin=*]
    \item  \textbf{Gset}: 67 binary-valued matrices, produced by an autonomous, random graph generator. The dimension $n$ varies from $800$ to $10^4$. \vspace{-0.2cm}
    \item \textbf{Dimacs10}: 152 symmetric matrices with $n$ varying from $39$ to $50,912,018$~, chosen for the $10$-th Dimacs Implementation Challenge.  We consider the 136 instances with dimension $n \le 8\cdot 10^6$. \vspace{-0.2cm}
\end{itemize}

\begin{figure}[t]
\centering
\vspace{-0.3cm}
\includegraphics[width=\columnwidth]{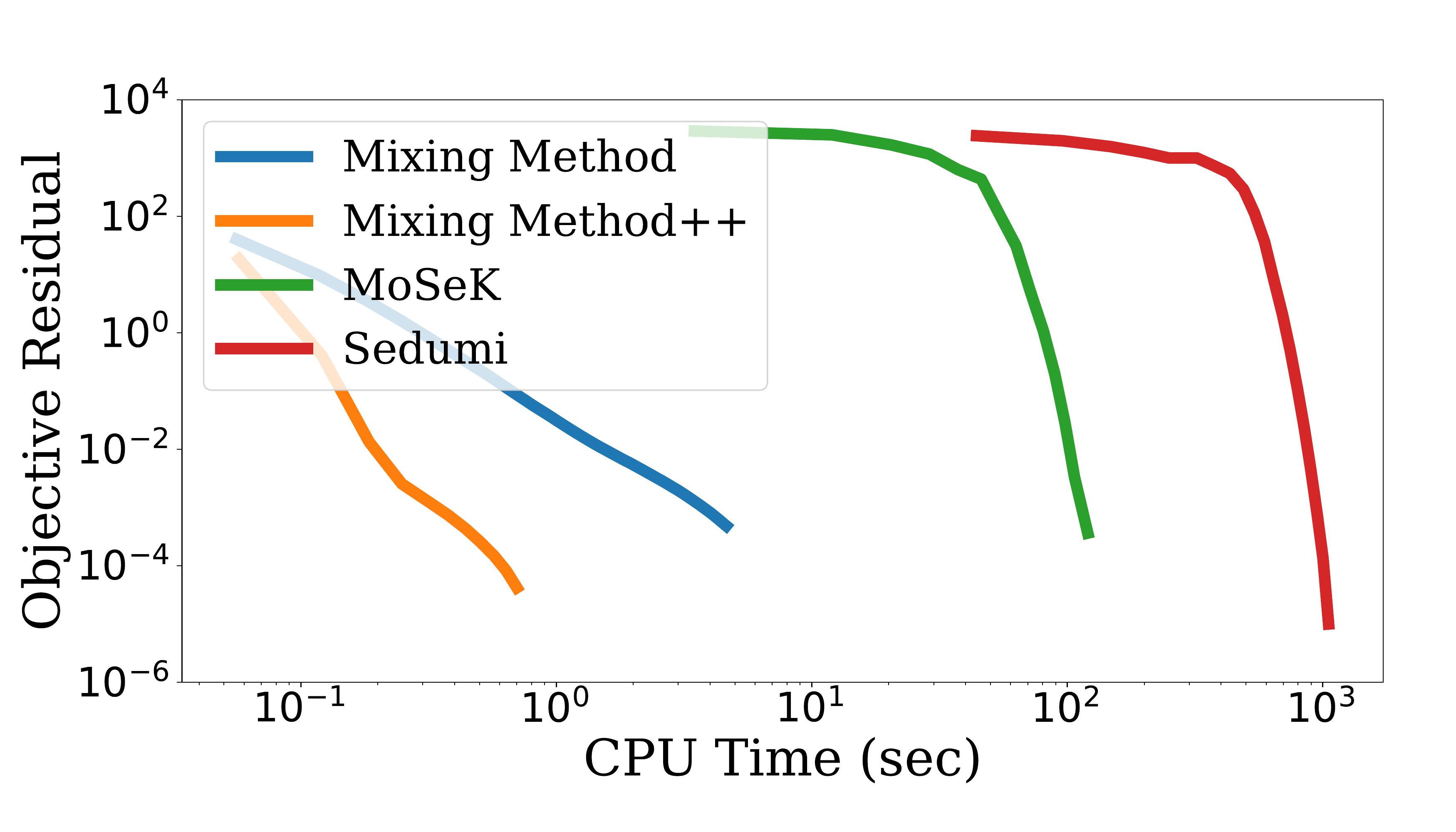}
\caption{\textbf{G40: $n=2000$, number of edges $ = 11766$.} \texttt{Mixing Method++} converges within $1$s while Mixing Method uses about $4$s. MoSeK and Sedumi need $10^2$ and $10^3$ seconds respectively to converge.} 
\label{exp:maxcut_G40}
\end{figure}

\begin{figure}[t]
\centering
\vspace{-0.3cm}
\includegraphics[width=\columnwidth]{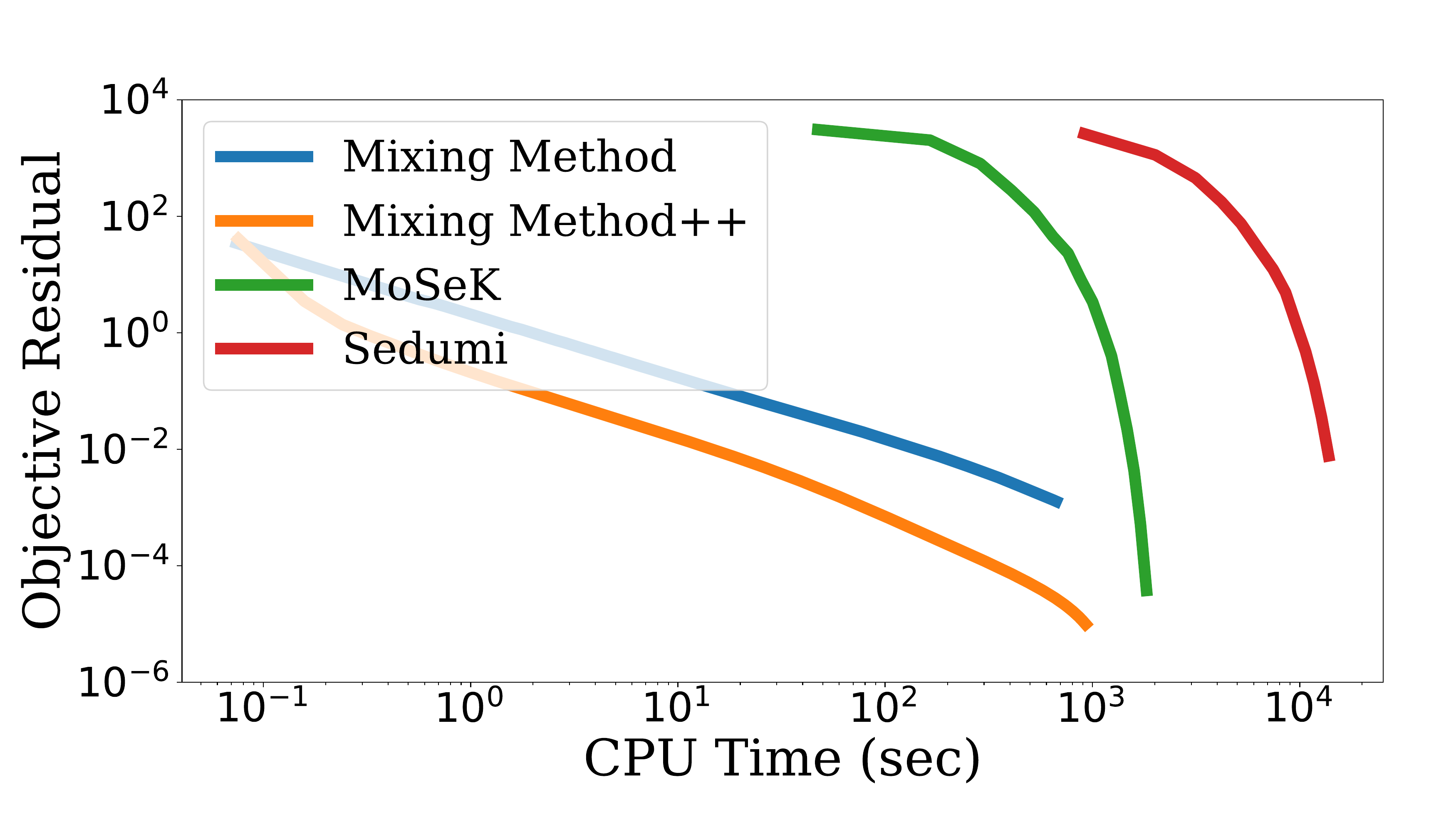}
\caption{\textbf{uk: $n=4824$, number of edges$ = 6837$.} Consider the CPU time when each solver reaches objective residual $10^{-2}$. \texttt{Mixing Method++} uses only 20s while Mixing Method needs 200s. MoSeK and Sedumi consume more than $10^3$ and $10^4$ seconds respectively.}
\label{exp:maxcut_uk}
\end{figure}

Table~\ref{table:maxcut} lists the number of solved instances, median logarithm of objective residual, and the pairwise acceleration ratio for different SDP solvers tested on MaxCut.\footnote{When an instance is solved by both \texttt{Mixing Method++} and solver $s$ within the 24 hour time limit, the pairwise acceleration ratio is defined as the CPU time taken by solver $s$ divided by that of \texttt{Mixing Method++}.} 
In comparison to the other solvers, \texttt{Mixing Method++} solves the most MaxCut instances (i.e., 111 of 203 total instances) within the imposed 24 hour time limit and provides a $5\times$ to $1316\times$ speedup.
Moreover, the efficiency of \texttt{Mixing Method++} does not harm its precision, which is evident in its median logarithm objective residual of $-3.87$.
Such performance is preceded only by Sedumi $(-4.30)$ and SDPNAL+ $(-4.33)$, both of which are more than $300\times$ slower than \texttt{Mixing Method++}.

We also study two specific cases of MaxCut: \textit{G40} from Gset and \textit{uk} from Dimacs10.
In Figures~\ref{exp:maxcut_G40} and~\ref{exp:maxcut_uk}, we plot objective residual with respect to CPU time of \texttt{Mixing Method++} and other SDP solvers on the \textit{G40} and \textit{uk} instances, respectively.\footnote{We do not include CGAL and SDPLR in the plots, because the intermediate solution given by those methods is infeasible, which makes intermediate objective values difficult to quantify.}
In both MaxCut instances, \texttt{Mixing Method++} converges significantly faster than other solvers that were considered. 

\subsection{MaxSAT}
\begin{table}
\centering
\begin{tabular}{ccc}
\toprule
\multicolumn{1}{c}{Methods} & & \multicolumn{1}{c}{Avg. Approx. Ratio} \\ 
\cmidrule(lr){1-1} \cmidrule(lr){3-3}
\multicolumn{1}{c}{Loandra} & & \multicolumn{1}{c}{0.945} \\ 
\multicolumn{1}{c}{Mixing Method} & & \multicolumn{1}{c}{0.975} \\ 
\multicolumn{1}{c}{\texttt{Mixing Method++}} & & \multicolumn{1}{c}{\textbf{0.977}} \\ \bottomrule
\end{tabular}
\caption{Results on MaxSAT. \texttt{Mixing Method++} achieves the best approximation ratio with a time limit of 300s. 850 problems are tested in total. Problems are derived from crafted (331), random (454), and industrial domains (65).}
\label{table:maxsat}
\end{table}


Let $s_{ij}\in\{\pm 1,0\}$ be the sign of variable $i$ in clause $j$.
The following problem provides an upper bound to the exact MaxSAT solution: \vspace{-0.05cm}
$$
\underset{\|v_i\|=1,}{\mathrm{max}} \sum_{j=1}^m\left(1-\tfrac{\|Vs_j\|^2-(|s_j|-1)^2}{4|s_j|}\right), \vspace{-0.05cm}
$$
where $s_j=[s_{1j},...,s_{nj}]^\top$.
Therefore, the MaxSAT SDP relaxation can then be solved by minimizing the following, related objective:
$$
f_{\text{MaxSAT}} = \left\langle C, {V^\top V} \right \rangle~~\text{s.t.}~~\|v_i\|_2=1~~\forall i,
$$
where $V \in \mathbb{R}^{k \times n}$ and $C=\sum_{j=1}^m\frac{s_js_j^\top}{4|s_j|}$.

We consider 850 problem instances from the 2016-2019 MaxSAT competition.
Each instance is categorized as \emph{random, crafted} or \emph{industrial}.
We evaluate \texttt{Mixing Method++} as a partial MaxSAT solver\footnote{A partial solver is only required to provide the assignment it finds with the least number of violated clauses, while a complete solver must also provide proof that such an assignment is optimal.} and compare its performance to that of the Mixing Method and Loandra; the latter was the best partial solver in the 2019 MaxSAT competition.
The best solution given by all solvers is used as ground truth when the optimal solution is not known.

\begin{figure}[t]
\centering
\includegraphics[width=\linewidth]{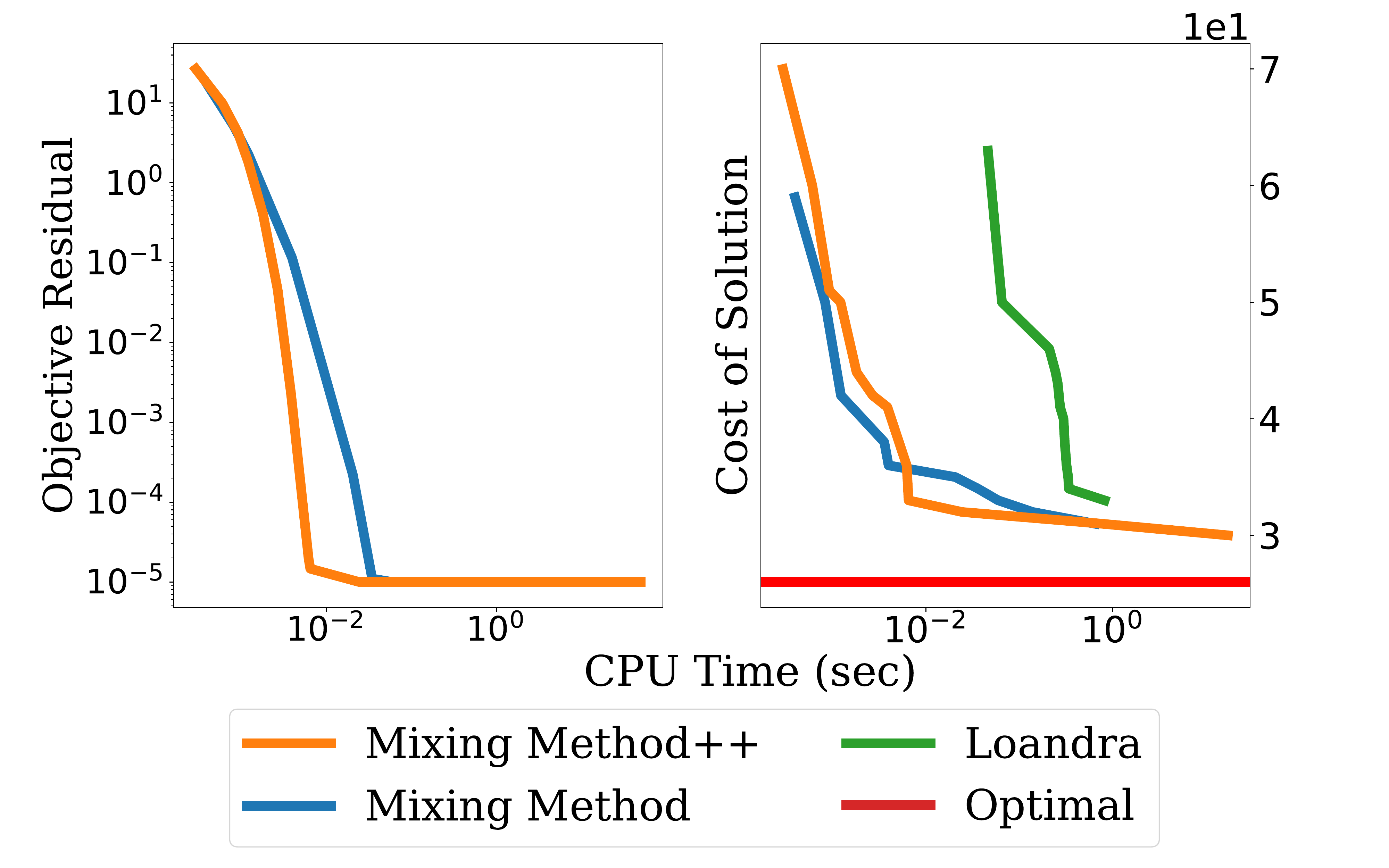}
\caption{\textbf{s3v90c800-6.cnf ($n=90$, $m=800$).} \texttt{Mixing Method++} not only converges faster than Mixing Method, but provides better solution (30) compared to Mixing Method (31) and Loandra (33). The optimal solution has cost 26.}
\label{exp:maxsat_s3v}
\end{figure}

\begin{figure}[t]
\centering
\includegraphics[width=\linewidth]{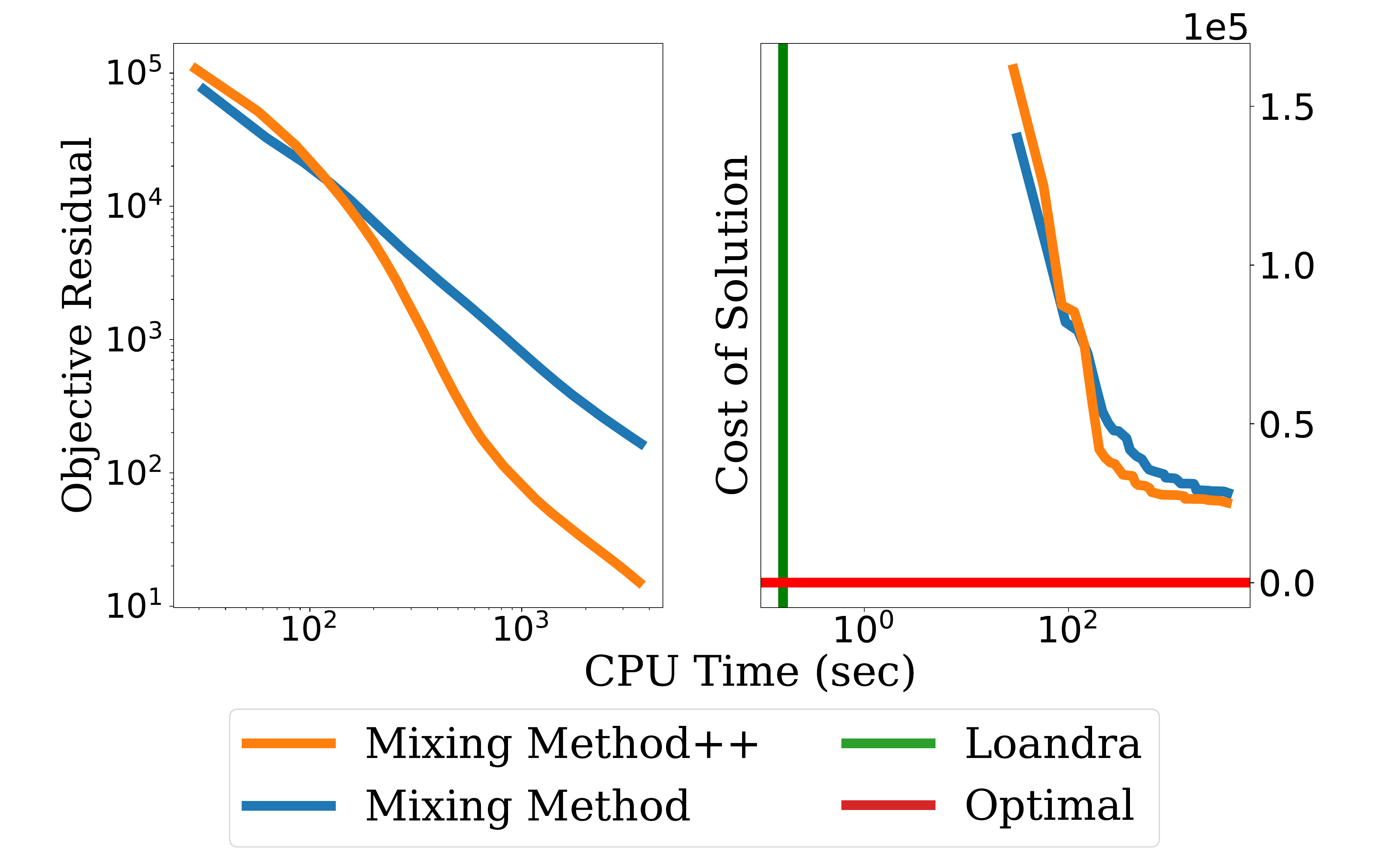}
\caption{\textbf{wb\_4m8s4\_dimacs\_filtered.cnf ($n=463,080$, $m=1,759,150$).} Loandra reaches optimal (230) in 1 sec, while Mixing Method and \texttt{Mixing Method++} return a solution with cost $28,306$ and $25,210$ respectively in 1 hour.}
\label{exp:maxsat_wb}
\end{figure}

The average approximation ratio of each solver across all MaxSAT instances is shown in Table~\ref{table:maxsat}.
\texttt{Mixing Method++} achieves the best average approximation ratio.
We also consider two specific MaxSat instances: \textit{s3v90c800-6} from the random track and from the industrial track \textit{wb.4m8s4.dimacs.filtered}.
In Figures~\ref{exp:maxsat_s3v} and~\ref{exp:maxsat_wb}, we plot both the objective residual and the cost (i.e., the number of unsatisfied clauses) with respect to CPU time for \textit{s3v90c800-6} and \textit{wb.4m8s4.dimacs.filtered} instances, respectively.
\texttt{Mixing Method++} converges faster than Mixing Method in both cases and returns the best solution on the \textit{s3v90c800-6} instance.
Although Loandra provides the best solution to \textit{wb.4m8s4.dimacs.filtered} quickly, \texttt{Mixing Method++} still outperforms Mixing Method significantly.


\subsection{MIMO}
The MIMO signal detection setting is defined as follows:
\begin{equation*}
        \argmin\limits_{x \in \mathcal{S}} \|y - Hx\|_2^2,
        \label{mimo_obj}
\end{equation*}
where $y \in \mathbb{R}^m$ is the signal vector, $H \in \mathbb{R}^{m \times n}$ is the channel matrix, and $x \in \mathcal{S}$ is the transmitted signal.
It is assumed that $y = Hx + v$, where $v$ is a noise vector.
For our experiments, we consider the constellation $\{ \pm 1\}$, which yields $\mathcal{S} = \{\pm1\}^n$.
The MIMO SDP relaxation can then be solved by minimizing the following objective.
\begin{align*}
f_{\text{MIMO}} = \langle C, V^TV \rangle~~\text{s.t.}~~\|v_i\|_2 = 1~~\forall i,
\end{align*}
where $V \in \mathbb{R}^{k \times (n + 1)}$ and $C = \begin{bmatrix} H^TH & -H^Ty \\ -y^TH & y^Ty\end{bmatrix}$.

\begin{figure}
    \centering
    \includegraphics[width=1\linewidth]{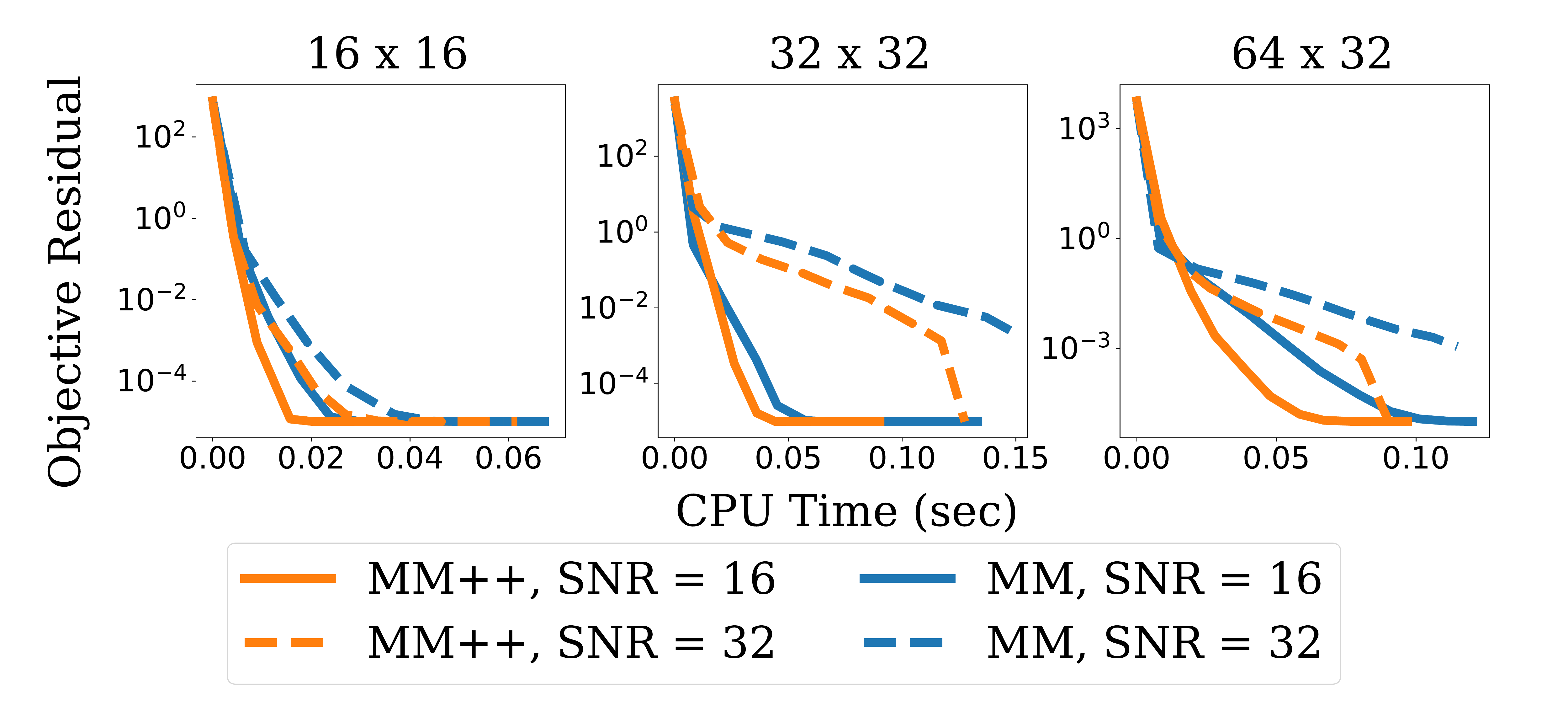}
    \caption{Convergence comparison of \texttt{Mixing Method++} and Mixing Method on MIMO signal detection problems with different SNRs. The size of the channel matrix is listed above each subplot.}
    \label{fig:mimo_mm}
\end{figure}

We compare \texttt{Mixing Method++} to Mixing Method on several simulated instances of the MIMO SDP objective.
$H$ is sampled from a standard normal distribution, $x$ is sampled from $\{\pm1\}$, and $v$ is sampled from a centered normal distribution with variance $\sigma_v^2$.
Experiments are performed for numerous settings of $\sigma_v^2$ determined by the signal-to-noise ratio (SNR) (i.e.,  $\sigma_v^2 = \frac{mn}{\text{SNR}}$).
We perform tests with $\text{SNR} \in \{8, 16\}$ and problem sizes $(m, n) \in \{(16, 16), (32, 32), (64, 32)\}$.
We plot the objective residual achieved by both Mixing Method and \texttt{Mixing Method++} on MIMO instances with respect to CPU time in Figure~\ref{fig:mimo_mm}.
As can be seen, \texttt{Mixing Method++} matches or exceeds the performance of the Mixing Method in all experimental settings.

\subsection{The Effect of Sparsity}
\begin{figure}[!ht]
    \centering
    \includegraphics[width=1\linewidth]{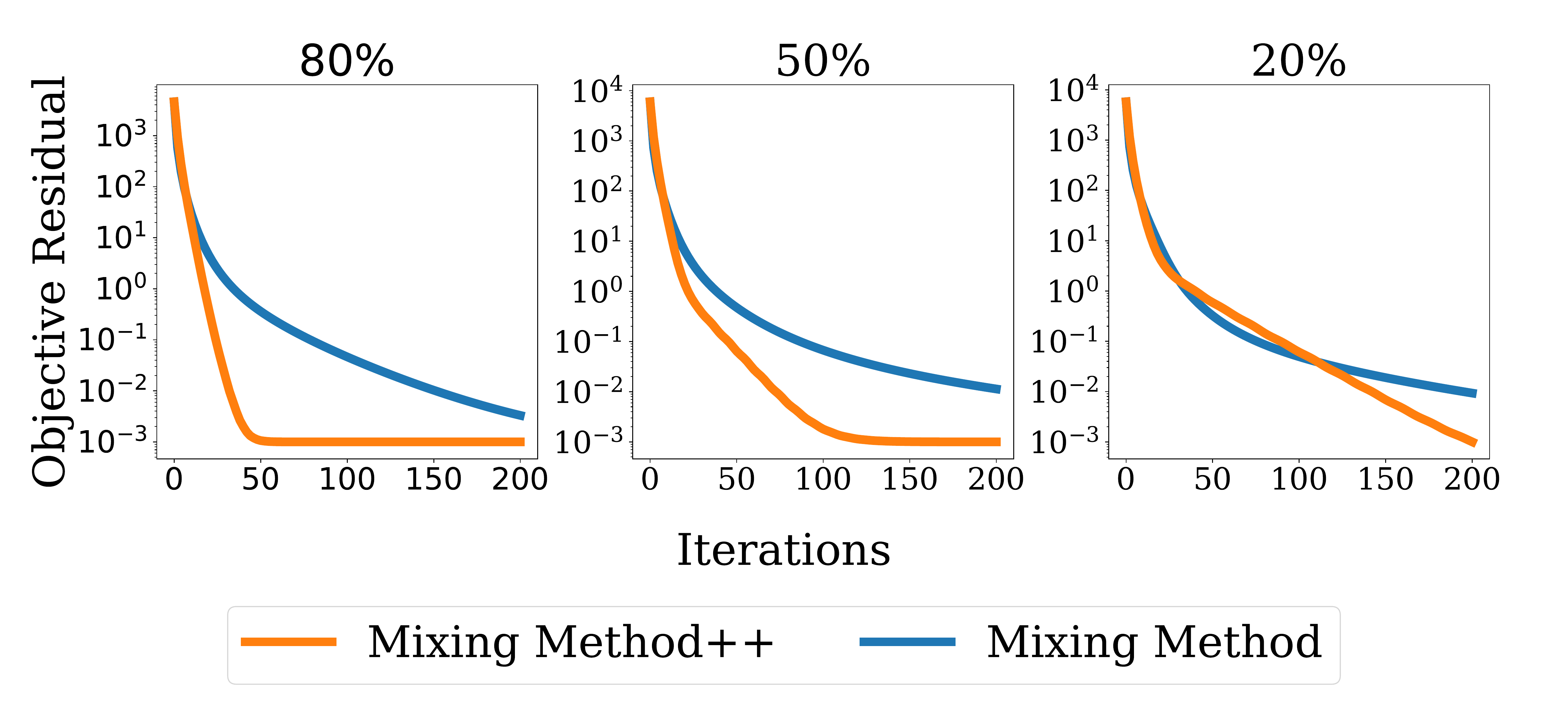}
    \caption{Convergence comparison of \texttt{Mixing Method++} (blue line) with $\beta=0.8$ and Mixing Method (orange line) on 500-node MaxCut problem instances with different sparsity levels.}
    \label{fig:sparsity}
\end{figure}

We empirically observe that \texttt{Mixing Method++} achieves different levels of acceleration based on the sparsity of $C$.\footnote{We define sparsity as the number of zero-valued elements divided by the total number of elements.}
To further understand the impact of sparsity on the performance of \texttt{Mixing Method++}, we compare \texttt{Mixing Method++} to Mixing Method for objectives with different sparsity levels.
In particular, we construct a MaxCut objective for a graph with 500 nodes and consider adjacency matrices with three different sparsity levels: $0.8$, $0.5$, and $0.2$.
In Figure~\ref{fig:sparsity}, we plot the objective residual of both Mixing Method and \texttt{Mixing Method++} on these problem instances.
As can be seen, the acceleration achieved by \texttt{Mixing Method++} is more significant when $C$ is sparse.
Furthermore, decreasing the value of $\beta$ improves the performance of \texttt{Mixing Method++} on MaxCut instances when $C$ is dense.

\subsection{Adaptive Momentum}
\begin{figure}[!ht]
    \centering
    \includegraphics[width=0.9\linewidth]{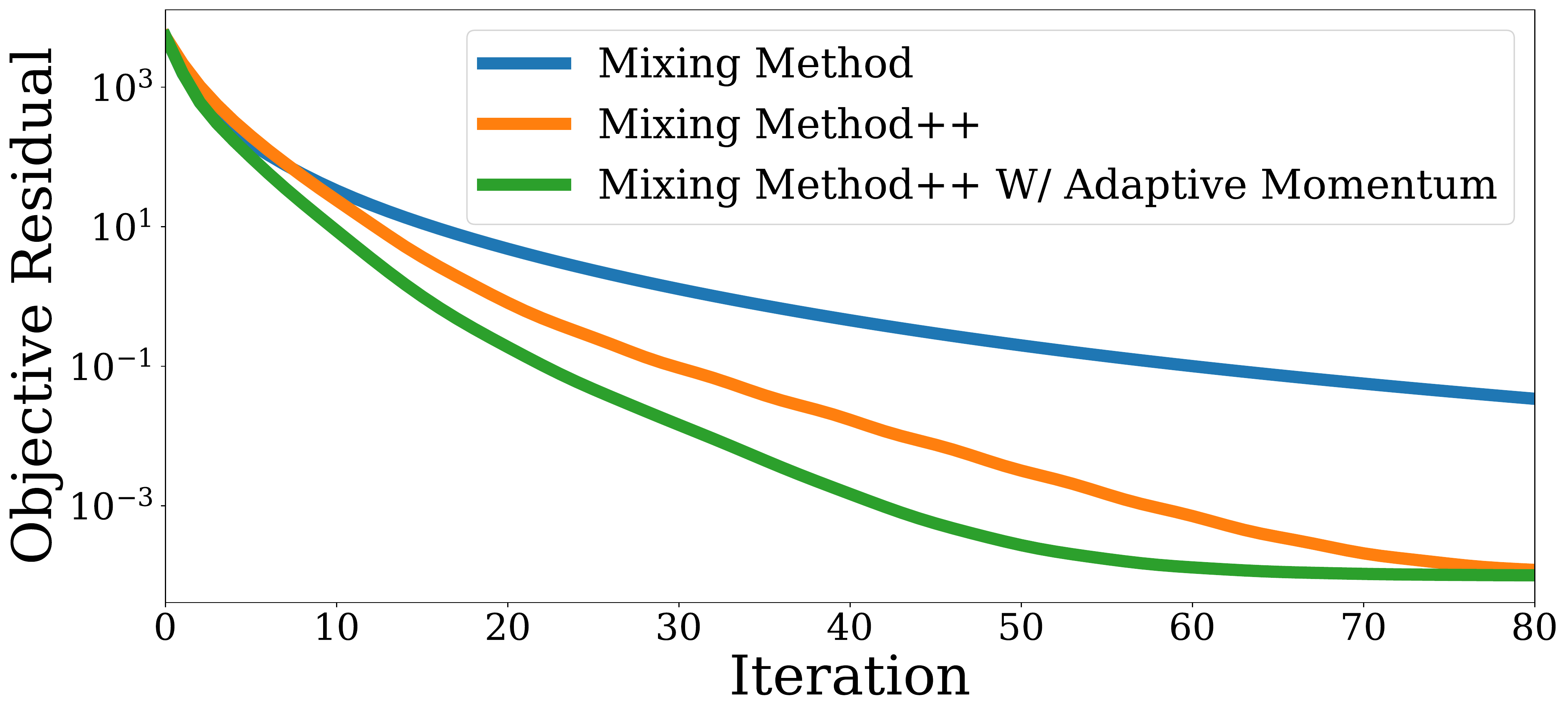}
    \caption{Comparison of Mixing Method (orange line), \texttt{Mixing Method++} with fixed momentum $\beta=0.8$ (blue line), \texttt{Mixing Method++} with an increasing Scheduled Momentum from $0$ to $0.8$ (green line) in terms of the objective cost per iteration.}
    \label{fig:scheduled-momentum}
\end{figure}

We observe that Mixing Method sometimes decreases the objective residual faster than \texttt{Mixing Method++} during early stage iterations.
Inspired by this observation and previous work on adaptive momentum schedules \cite{o2015adaptive,chen2019demon,wang2020scheduled}, we explore whether dynamically adjusting the value of $\beta$ between iterations can further improve the performance of \texttt{Mixing Method++}.
We propose the following momentum schedule:
\begin{align*}
    \beta_t = \beta \left(1-\exp(-\alpha t/T)\right),
\end{align*}
where $\beta = 0.8$ is the default momentum, $t$ is the current iteration, $T$ is the total number of iterations, and $\alpha$ is a constant such that $\exp(-\alpha) \approx 0$. 
Intuitively, this schedule begins with a low momentum value, which quickly increases to $\beta=0.8$, thus forming a warm-up schedule for the momentum.
As shown in Figure \ref{fig:scheduled-momentum}, this schedule improves the convergence speed of \texttt{Mixing Method++}.


\section{Local linear convergence}\label{sec:proofs}

We now give detailed proof of Theorem~\ref{theorem: linearconvergence}.
Let $g_i$ be the linear combination in step~\ref{step:mixing inner update} of Algorithm~\ref{algo:mixing}, before projecting onto the unit sphere:
\begin{equation*}\label{eq: g_i}
\begin{split}
    g_i &= \sum_{j< i} c_{i j} \hat{v}_j  + \sum_{j> i} c_{i j} v_j.
\end{split}
\end{equation*}
We also define $y_i = \|g_i\|_2$ and $f(V) \triangleq \langle C, V^TV \rangle$ as the objective value for some viable solution $V$.
Therefore, step~\ref{step:mixing inner update} of Algorithm~\ref{algo:mixing} can be expressed as: $y_i u_i = - g_i$. 
Thus, the changes presented in the Algorithm~\ref{algo:accelerated_mixing} can be written as the combination of~\ref{step:mixing++ inner update} and the momentum, in step~\ref{step:mixing++ inner momentum}:
\begin{equation*}\label{eq:twoSteps}
\begin{split}
    u_i = - \frac{g_i}{y_i}, \hspace{10 pt} \text{ and } \hspace{10 pt} 
    \hat{v}_i = \frac{\left(1+\beta\right)u_i - \beta v_i}{w_i},
\end{split}
\end{equation*}
where \mbox{$w_i = \left\lVert \left(1+\beta\right)u_i - \beta v_i\right\rVert_2$}. 
This leads to the following one line relation,
\begin{equation}\label{eq:recursion}
    w_iy_i\hat{v}_i  = -\left(1+\beta\right)g_i -\beta y_i v_i, ~\forall i \in [n].
\end{equation}

Herein, the matrices $V$ and $\hat{V}$ refer to the current and the end state of the inner iteration in algorithm~\ref{algo:accelerated_mixing}. For each $i\in [n]$, we define the family of maps \mbox{$Z_i : (S^{k-1})^{n} \to (S^{k-1})^{n}$} as follows
\begin{equation}\label{eq:Z}
    Z_i(V)\, e_j = \begin{cases}
\hat{v}_j, \hspace{10 pt} \text{if $j< i$}\\
v_j, \hspace{10 pt} \text{if $j\ge  i$},
\end{cases}
\end{equation}
where $\{e_j\}$ is the canonical basis in $\mathbb{R}^n$. It is trivial to see that $Z_i(V)$ corresponds to the matrix in the inner cycle before updating $v_i$ to $\hat{v}_i$. Thus, only the vector $v_i$ has been changed from $Z_i(V)$ to \mbox{$Z_{i+1}(V)$}. Also notice that \mbox{$Z_{n+1}(V)$} corresponds to the end of the inner cycle, i.e., \mbox{$Z_{n+1}(V) = \hat{V}$}.

\subsection*{Matrix Representation}
The transition from the state $v_i$ to $\hat{v}_i$ is expressed in the relation~\eqref{eq:recursion}, then for any $i \in [n]$ in each inner cycle we have, 
\begin{align}\label{eq:matrixA2}
\begin{split}
   \left(1+\beta\right)\sum_{j< i} c_{i j} \hat{v}_j + w_i y_i \hat{v}_i =-\left(1+\beta\right)\sum_{j> i} c_{i j}  v_j - \beta y_i v_i.
\end{split}
\end{align}
Let $L$ be the strictly lower triangular matrix of $C$ so~\eqref{eq:matrixA2} leads us to the representation: 
\begin{align}\label{eq:matrixformA2}
    \left(\left(1+\beta\right)L+ D_w D_y\right)\, \hat{V}^\top = - \left(\left(1+\beta\right)L^\top + \beta D_y\right)\, V^\top.
\end{align}

\subsection*{Proof of Lemma~\ref{lemma:decreasing}}
Let $\beta$ be in $[0, 1]$.  When fixing all the others variables, $f$ with respect to $v_i$ is given by
\begin{equation*}
     \langle C,V^\top V \rangle= 2 v_i^\top \big(\sum_j c_{ij}v_j \big)+ \text{constant}.
\end{equation*}
Since only $v_i$ is changed, during the $i$th coordinate-wise updating, the only part of the objective that will change is $2v_i^\top g_i$. That is,
\begin{equation*}
 f(Z_i(V))-f(Z_{i+1}(V)) = 2g_i^\top \left(v_i-\hat{v}_i\right). 
\end{equation*}

The updating rule of algorithm~\ref{algo:accelerated_mixing} is given by the recursion~\eqref{eq:recursion}, so using this into the above equation yields
\begin{equation*}
\begin{split}
   f(Z_i(V))-f(Z_{i+1}(V)) &= \frac{2y_i\,(w_i- \beta)}{1+\beta}(1-v_i^\top \hat{v}_i),\\
    &= \frac{y_i\, (w_i- \beta)}{1+\beta}\|v_i - \hat{v}_i\|_2^2.
\end{split}
\end{equation*}

From observation~\ref{observation0}, we have  $w_i-\beta \ge 1-\beta$. Thus,
\begin{equation*}
\begin{split}
    f(Z_i(V))-f(Z_{i+1}(V)) \ge y_i\, \frac{1- \beta}{1+\beta}\|v_i - \hat{v}_i\|_2^2.
\end{split}
\end{equation*}
Adding over all the $i$, leads us to the desired conclusion 
\begin{equation*}
\begin{split}
    f(V)-f(\hat{V}) &= \sum_{i= 1}^n\,f(Z_i(V))-f(Z_{i+1}(V)) \\
    &\ge \frac{1-\beta}{1+\beta}\sum_{i= 1}^n y_i\,\|v_i - \hat{v}_i\|_2^2.
    \end{split}
\end{equation*}

The following observation will be useful later on.

\begin{observation}\label{observation}
For all $i\in[n]$, $w_i^*=1$
\end{observation}
\begin{proof} Since $v^*_i$ is a fixed point, using this fact in Eq.~\eqref{eq:recursion} yields \mbox{$(w^*_i+\beta)v_i^*=u_i^*$}. Taking the norm, \mbox{$w_i^*+\beta=1+\beta$}, since $u^*_i$ and $v^*_i$ are in the unit sphere.
\end{proof}

Let $S$ be defined as $\left(1+\beta\right)C+D_wD_y + \beta D_y$ and let $S^*$ be $\left(1+\beta\right)C+D_{w^*}D_{y^*} + \beta D_{y^*}$; namely, the corresponding value of $S$ when $V$ is optimal in the optimization~\eqref{eq:noncvx_continuous_form}. We have the following observation.
\begin{observation}\label{PSD}  $S^*$ is PSD.
\end{observation}

\begin{proof}
By observation~\ref{observation} we have $w^*_i = 1$, so 
\begin{equation*}
  S^* = \left(1+\beta\right)C+D_{w^*}D_{y^*} + \beta D_{y^*} = \left(1+\beta\right)\left(C+ D_y^*\right),
\end{equation*}
where according to lemma~3.12, and since $V$ is optimal, it follows that $\frac{S^*}{1+\beta}\succeq 0$, hence $S^*\succeq 0$. 
\end{proof}

\subsection*{Proof of Lemma~\ref{lemma: lemmad2}}

By~\eqref{eq:matrixformA2}, we have
\begin{align*}
    \hat{V} =\,- V\left(\left(1+\beta\right)L + \beta D_y\right) \left(\left(1+\beta\right)L^\top +D_w D_y\right)^{-1}.
\end{align*}
This leads to
\begin{align*}
\begin{split}
 V -\hat{V}
 =& V(\left(1+\beta\right)C + D_wD_y + \beta D_y)\\
 &(\left(1+\beta\right)L^\top + D_wD_y) ^{-1}\\
 =& V S \underbrace{(\left(1+\beta\right)L^\top + D_wD_y) ^{-1}}_{R},
     \end{split}
\end{align*}

Using the definition \mbox{$\Delta\triangleq (D_wD_y+\beta D_y)- (D_{w^*}D_{y^*}+\beta D_{y^*})$}, and the fact that \mbox{$S =  S^*+\Delta$} in the previous equation, we get
\begin{equation}\label{eq:inequalityAlg2}
    \begin{split}
        \|V-\hat{V}\|_F^2 
        &\ge  \|VS^*R\|_F^2 + 2 \text{tr}\big(V^\top VS^* R R^\top \Delta\big).
    \end{split}
\end{equation}
We claim that $$\|(w+\beta \mathds{1})\odot y- (w^*+\beta \mathds{1})\odot  y^*\|_2 I_n + \Delta	\succeq 0.$$
To see this, first notice \mbox{$D_wD_y+ \beta D_y = D_{(w+\beta \mathds{1})\odot y}$}. Thus,
$$ \Delta = D_{(w+\beta \mathds{1})\odot y}- D_{(w^*+\beta \mathds{1})\odot y^*}.$$
Besides, 
$$\|(w+\beta \mathds{1})\odot y- (w^*+\beta \mathds{1})\odot  y^*\|_2 = \| \Delta \|_F.$$
So the claim follows since $\beta \ge 0$ and termwise 
\begin{equation*}
    \begin{split}
        \| \Delta \|_F+ \Delta_i \ge |\Delta_i |+ \Delta_i \ge 0. 
    \end{split}
\end{equation*}
Now, for the first part in~\eqref{eq:inequalityAlg2}, we get 
\begin{equation*}
\begin{split}
    \|VS^*R\|_F^2 &= \text{tr}(\underbrace{{S^*}^\top}_{\succeq 0} \underbrace{V^\top V}_{\succeq 0} \underbrace{S^*}_{\succeq 0} \underbrace{R R^\top}_{\succeq 0}) \\
    &\ge \sigma_{\min}^2(R) \sigma_{\texttt{nnz}}(S^*) \text{tr}(V^TV S^*).
\end{split}
\end{equation*}
For the second in~\eqref{eq:inequalityAlg2}, we get
\begin{equation*}
    \begin{split}  \text{tr}\big(V^\top VS^* R R^\top \Delta\big) &\ge -\|\Delta\|_F\tr\big(V^\top VS^* R R^\top\big)\\
    &\ge-\|\Delta\|_F\,\sigma_{\max}^2(R)\text{tr}\big(V^\top VS^*\big),
    \end{split}
\end{equation*}
this last inequality is given by our previous claim, the fact that $S^* \succeq 0$ (observation~\ref{PSD}) and 
\begin{equation*}
    \begin{split}
        \tr(V^\top VS^*RR^\top) &=\tr(V^\top VS^*\underbrace{UD^2U^\top}_{RR^\top}) =\tr(V^\top VS^*D^2)\\
        &\le \sigma^2_{\max}(R) \, \tr(V^\top VS^*).
    \end{split}
\end{equation*}
Using the lower bounds in~\eqref{eq:inequalityAlg2} gives us
\begin{equation*}
    \begin{split}
        \|V&-\hat{V}\|_F^2 \ge \big(\sigma_{\min}^2(R) \sigma_{\texttt{nnz}}(S^*) -2\|\Delta\|_F\sigma_{\max}^2(R)\big)\text{tr}\big(V^\top VS^*\big).
    \end{split}
\end{equation*}
To conclude, the term 
\begin{equation*}
    \begin{split}
      \tr(V^\top VS^*) &= \tr(V^\top VC)-\tr(-V^\top V D_{y^*} ) 
      = f(V) - \tr(-V^\top V D_{y^*})\\
      &\stackrel{i)}{=} f(V) + \sum\limits_{i} \underbrace{\langle v_i, v_i \rangle}_{=1} \cdot y^*_i
      = f(V) + \langle \mathds{1}, y^* \rangle\\   
      &\stackrel{ii)}{\geq} f(V) - f^*.
    \end{split}
\end{equation*}
where $i)$ follows from the fact that $D_{y^*}$ is diagonal. 
To understand $ii)$, consider the dual to the diagonally constrained SDP problem~\eqref{eq:continuous_form} (i.e., see Lemma 3.12 of \cite{wang2017mixing} for more details):
\begin{align}
    \max\limits_{y} -\langle \mathds{1}, y \rangle, ~\text{such that } C + D_y \succeq 0
    \label{dual_sdp}
\end{align}
We know that \mbox{$C + D_{y^*} = \frac{S^*}{1+\beta} \succeq 0$} from Observation~\ref{PSD}.
Therefore, $y^*$ is a feasible solution to the dual problem~\eqref{dual_sdp}, revealing that \mbox{$f^* \geq - \langle \mathds{1}, y^* \rangle$}.
This fact allows us to derive the inequality given by $ii)$.
Then, if we define \mbox{$\zeta = \sigma^2_{\min} (R)\sigma_{nnz} (S^*)$} and \mbox{$\gamma = 2\sigma^2_{\max} (R)$} with \mbox{$\sigma_{\min}(R) = 1/(w\odot y)_{\max}$} and \mbox{$\sigma_{\max}(R) = 1/(w\odot y)_{\min}$}, we arrive at the final result:
\begin{align*}
    \|V - \hat{V}\|^2_F \geq \left(\zeta - \gamma \|\Delta\|_F\right)(f(V) - f^\star)
\end{align*}
It should be noted that both $\zeta$ and $\gamma$ are strictly positive constants because $R$ is defined as an invertible matrix and, therefore, has all non-zero eigenvalues.

\subsection*{Proof of Lemma~\ref{lemma: lemmad4}}
Recall that $\Delta_i = (w_i + \beta) y_i - (w_i^\star + \beta) y_i^\star$.
We can unroll the value of $\Delta_i \hat{v}_i$ into the following expression: 
\begin{align}\label{eq:difference}
        \big((w_iy_i -w_i^*y_i^*)+ \beta(y_i-y_i^*)\big)\, \hat{v_i}
\end{align}
Then, we derive the following equality from~\eqref{eq:recursion}:
\begin{equation}
\begin{split}
    w_i y_i\hat{v_i} &= -\left(1+\beta\right)g_i - \beta y_i v_i \\
    &= -\left(1+\beta\right)Z_i c_i - \beta y_i v_i \label{eq:recursion_unrolled}
\end{split}
\end{equation}
where we denote \mbox{$Z_i = Z_i(V)$} for brevity and $c_i$ as the $i$th column of $C$.
We denote \mbox{$s_i^*  = \left(1+\beta\right)c_i +w_i^*y_i ^* e_i + \beta y_i ^* e_i$}.
By adding~\eqref{eq:difference} to both sides of~\eqref{eq:recursion_unrolled}, we arrive at the following, expanded version of the update rule in~\eqref{eq:recursion}.
\begin{equation}\label{eq:eqSubs}
    \begin{split}
        &w_iy_i \hat{v_i} -w_i^*y_i^* \hat{v_i}+ \beta y_i \hat{v_i} - \beta y_i^*\hat{v_i} \\
        &= -\left(1+\beta\right)Z_i c_i - \beta y_i v_i -w_i^*y_i^* \hat{v_i}+ \beta y_i \hat{v_i} - \beta y_i^*\hat{v_i}\\
        &\stackrel{i)}{=} -Z_i s_i^ *+ \big(w_i^*y_i^* - \beta (y_i -y_i^ *) \big) (v_i - \hat{v_i}).
    \end{split}
\end{equation}
where $i)$ follows from the definition of $s_i^\star$ and the fact that $Z_i e_i = v_i$.
By combining Assumption~\ref{assu:delta} and Lemma~\ref{lemma:decreasing}, it is known that
\begin{equation}\label{eq:ineq1}
    f(V)- f^ * \ge   \frac{1-\beta}{1+\beta}\,\delta \, \|v_i - \hat{v_i} \|_2^2.
\end{equation}
Besides,
\begin{align}\label{eq:ineq2}
\begin{split}
    f(V)- f^* &\ge f(Z_i)- f^*
    = \tr(Z_i^\top Z_i S^*) \\ &\ge \sigma^{-1}_{\max}(S^*)\|Z_i S^*\|_F^2
    \ge \sigma^{-1}_{\max}(S^*)\|Z_i s_i^*\|_2^2
\end{split}
\end{align}

Therefore, based on~\eqref{eq:ineq1} and~\eqref{eq:ineq2}, it is known that \mbox{$\|Z_i s_i^\star\|_2^2$} and \mbox{$\|v_i - \hat{v}_i\|^2_2$} are both upper bounded by a positive factor of \mbox{$f(V) - f^\star$}.
The inner product of~\eqref{eq:difference} with itself to yields: 
\begin{align*}
    \Delta_i^2 &\stackrel{i)}{=} \big((w_iy_i -w_i^*y_i^*)+ \beta(y_i-y_i^*)\big)^2\\
    &\stackrel{ii)}{\le}2 \|Z_i s_i^ *\|^ 2  + 2 (y^*_i-\beta y_i)^2 \|v_i - \hat{v_i}\|^2 \\
    &=2 \|Z_i s_i^ *\|^ 2  + 2 ({y^*_i}^2-2\beta y_i^*y_i+ {y_i}^ 2)\, \|v_i - \hat{v_i}\|^2
\end{align*}
where $i)$ follows from the fact that $\hat{v}_i$ is a unit vector and $ii)$ follows from combining observation~\ref{observation} with~\eqref{eq:eqSubs}. 
Next, we have: $$\|(w+\beta \mathds{1})\odot y- (w^*+\beta \mathds{1})\odot  y^*\|_2^2 = \|\Delta\|_F^2 = \sum\limits_{i=1}^{n} \Delta_i^2.$$
We combine this expression with the upper bounds for \mbox{$\|Z_i s_i^\star\|^2_2$} and \mbox{$\|v_i - \hat{v}_i\|^2_2$} derived in~\eqref{eq:ineq1} and~\eqref{eq:ineq2} to arrive at the final result, where $\tau$ is defined as some positive constant: 
\begin{align*}
\|(w+\beta \mathds{1})\odot y- (w^*+\beta \mathds{1})\odot  y^*\|_2^2 \le \tau (f(V) - f^*)
\end{align*}

\subsection*{Proof of Theorem~\ref{theorem: linearconvergence}}
Based on Lemma~\ref{lemma: lemmad4}, a neighborhood can be selected around the optimum such that the value of $f(V) - f^\star$ is sufficiently small. 
We define this neighborhood through the selection of a positive constant $\kappa$ such that the following inequality holds, where $\tau$ is defined in Lemma~\ref{lemma: lemmad4}:
\begin{align}
    \left ( \frac{\zeta - \kappa}{\gamma}\right )^2 \geq \tau \left ( f(V) - f^\star \right) \label{neighborhood_exp}
\end{align}
Within this expression, $\zeta$ and $\gamma$ are both defined in Lemma~\ref{lemma: lemmad2}.
Because $\zeta$ and $\gamma$ are both strictly positive (i.e., see Lemma~\ref{lemma: lemmad2}), $\left ( \frac{\zeta - \kappa}{\gamma} \right )$ is known to be strictly positive so long as $\kappa \in [0, \zeta)$.
Therefore, there always exists a value of $\kappa$ such that~\eqref{neighborhood_exp} will be true within a sufficiently small neighborhood around the optimum. 
We combine~\eqref{neighborhood_exp} with the inequality from Lemma~\ref{lemma: lemmad4} to yield the following:
\begin{equation*}
    \left(\frac{\zeta - \kappa}{\gamma} \right)^ 2\ge \tau (f(V)- f^*) \ge \|\Delta\|_F^ 2,
\end{equation*}
Then, because $\left ( \frac{\zeta - \kappa}{\gamma} \right )$ is strictly positive, we can manipulate the above expression to yield the following:
\begin{equation}
    \zeta - \gamma\|\Delta\|_F  \ge \kappa. \label{kappa_bound}
\end{equation}

The above inequality holds for all successive iterations of Algorithm~\ref{algo:accelerated_mixing}.
Then, based on the combination of Lemma~\ref{lemma:decreasing} and assumption~\ref{assu:delta}, $f(V)- f(\hat{V}) \ge  \delta\frac{1-\beta}{1+\beta}\|V - \hat{V}\|_F^ 2$, and thus 
\begin{equation*}
\begin{split}
    f(V)- f(\hat{V}) &\ge  \delta \frac{1-\beta}{1+\beta}\|V - \hat{V}\|_F^ 2\\
                  &\stackrel{i)}{\ge}  \delta\frac{1-\beta}{1+\beta} \kappa (f(V)-f^*)
\end{split}
\end{equation*}
where $i)$ is derived by combining~\eqref{kappa_bound} and Lemma~\ref{lemma: lemmad2}.
This expression in turn implies the following.
\begin{equation*}
    \left(1- \delta \kappa \frac{1-\beta}{1+\beta}\right)(f(V)-f^*) \ge f(\hat{V})-f^*.
\end{equation*}
thus giving us the desired linear convergence.
By substituting $\rho = \delta \kappa \frac{1-\beta}{1+\beta}$, we arrive at the final result.

\section{Convergence to first-order critical point} \label{sec:critical}
\texttt{Mixing method++} not only has local linear convergence, but it always converges to a first-order critical point, per following Theorem: 

\begin{theorem}
Let $V^\ell = Z_{n+1}^\ell(V)$, see definition~\eqref{eq:Z}, for all $\ell\in \mathbb{N}$; i.e., the collection of points generated by \texttt{Mixing Method++} after finishing the inner loop. Then, under the assumption~\ref{assu:delta}, $\{V^\ell\}_\ell$ converges to a limit point $\bar{V}$ and the Riemannian gradient at $\bar{V}$ is zero (\texttt{Mixing Method++} converges to a first-order critical point).
\end{theorem}

\begin{proof}
For $\beta$ in $[0, 1)$ we know by lemma~\ref{lemma:decreasing} that the objective function is decreasing. Besides, $(S^{k-1})^n$ is compact, indeed limit point compact, so there exists a limit point $\bar{V}$ such that $f(\bar{V}) = \lim_k f(V^k)$ by continuity of $f$. It is clear that $\bar V$ is a fixed point in relation~\eqref{eq:matrixformA2}, so

\begin{equation}\label{eq:equilibrium}
   \bar{V} C = -\frac{\bar{V} D_{(\bar{w}+ \beta \mathds{1}) \odot \bar{y}}}{1+\beta}.
\end{equation}

Since the constraint set is a product of spheres, its corresponding tangent space (henceforth denoted by the letter $\mathrm{T}$) at the point $\bar{V}$ is given by the product of tangent spaces of $S^{k-1}$; namely

\begin{equation}\label{eq:ManProd}
    \mathrm{T}_{\bar{V}} \left(S^{k-1}\right )^n = \prod_{j=1}^n \mathrm{T}_{\bar{v_j}}\, \left(S^{k-1}\right),
\end{equation}
where $\bar{v}_j$ is the $j$th column of $\bar V$. 

The above result gives us a way to characterize the tangent space. To this end, it is actually enough to work with the tangent space of $S^{k-1}$. 

By definition, for any $v_j$, the vector $\dot{v}_j\in \mathbb{R}^k$ is in the tangent space of $S^{k-1}$ if and only if there exists a curve $\varphi: I \to S^{k-1}$ such that $\varphi(0) = v_j $ and $\dot{\varphi}(0) = \dot{v}_j$, where $\dot{\varphi}$ is the derivative of $\varphi$. Then $\varphi(t) \in S^{k-1}$ if satisfies $ \langle \varphi(t), \varphi(t) \rangle  = 1$ for all $t\in I$. Differentiating on both sides leads us to $0 = \langle  \dot{\varphi}(t), \varphi(t) \rangle  + \langle \varphi(t), \dot{\varphi}(t) \rangle$. Evaluating at $t = 0$, we get 
\begin{equation*}
     2 \langle  \dot{v}_j,  v_j \rangle  = 0,
\end{equation*}
thus $\{\dot{v}_j \in \mathbb{R}^k : \langle  \dot{v}_j,  v_j \rangle  = 0 \} \subseteq \mathrm{T}_{v_j} S^{k-1}$, since both subspaces are of the same dimension, $\{\dot{v}_j \in \mathbb{R}^k : \langle  \dot{v}_j,  v_j \rangle = 0  \}  = \mathrm{T}_{v_j} S^{k-1}$.

Now by~\eqref{eq:ManProd} and the previous analysis, we finally arrive to the characterization of the tangent space; namely 
\begin{equation*}
    \mathrm{T}_{\bar{V}} \left(S^{k-1}\right)^n = \left\{ \dot{\bar{V}} \in \mathbb{R}^{k\times n} :  \langle  \dot{\bar{v}}_j,  \bar{v}_j \rangle = 0 \hspace{5 pt} \forall j \in [n]\, \right\},
\end{equation*}
and $\dot{\bar{v}}_j$ ($v_j$) corresponds to the $j$th column of $\dot{\bar{V}}$ ($\bar{V}$), respectively.

Let $\mathrm{P}^\perp: \mathbb{R}^{k\times n}  \to  \mathrm{T}_{\bar{V}} \left(S^{k-1}\right)^n$ be the projector operator from the euclidean space to the tangent space at $\bar{V}$, defined for any $W \in \mathbb{R}^{k\times n}$ as follows,
\begin{equation*}
    W  \xrightarrow{ \mathrm{P}^\perp}  \begin{pmatrix}

   w_1 -   \langle  \bar{v}_1,  w_1 \rangle \bar{v}_1     & \ldots  & w_n -\langle  \bar{v}_n,  w_n \rangle \bar{v}_n
\end{pmatrix}.
\end{equation*}

Let the Riemannian gradient of $f$ denoted by $\mathrm{grad}f(\cdot)$, and the gradient of $f$ defined in the entire euclidean domain as $\nabla f (\cdot)$. Using tools from matrix manifold the Riemannian gradient is given by 

$$\mathrm{grad}f(\bar{V}) = \mathrm{P}^\perp (\nabla f(\bar{V})).$$  
It is easy to see that $\nabla f(\bar{V}) = 2 \bar{V}C$, and by relation~\eqref{eq:equilibrium}, this yields us to the equivalent version $\nabla f(\bar{V}) =  -\frac{2}{1+\beta} \bar{V} D_{(\bar{w}+ \beta \mathds{1}) \odot \bar{y}}$. 

Finally, \mbox{$\left(\nabla f(\bar{V})\right)_j = -\frac{2}{1+\beta}(w_j y_j+ \beta y_j) \bar{v}_j$}, \mbox{$\forall j\in [n]$} and by definition of the projection map we get 
\begin{equation*}
\begin{split}
     \left( \mathrm{P}^\perp(\nabla f(\bar{V})) \right) _j  =  \frac{2}{1+\beta}(w_j y_j + \beta y_j  )  (- \bar{v}_j  +\bar{v}_j) = 0.
\end{split}
\end{equation*}
So, we have $\left( \mathrm{P}^\perp(\nabla f(\bar{V}) \right) _j = 0$, for all coordinates $j\in[n]$, and thus
\begin{align*}
    \mathrm{grad}f(\bar{V}) = \mathrm{P}^\perp (\nabla f(\bar{V})) = 0.
\end{align*}
\end{proof}

\section{Conclusion and Future Directions} 
We present a novel approach, \texttt{Mixing Method++}, to solve  diagonally constrained SDPs. 
\texttt{Mixing Method++} inherits the simplicity---and non-trivially preserves theoretical guarantees from---its predecessor. 
In practice, it yields not only faster convergence on nearly all tested instances, but also improvements in solution quality.
\texttt{Mixing Method++} adds one addition hyperparameter $\beta$ for which we provide a theoretical upper bound (i.e., $\beta < 1$).
Using $\beta = 0.8$ experimentally leads to a robust algorithm, which outperforms Mixing Method and numerous other state-of-the-art SDP solvers.
Theoretically validating the experimental acceleration provided by \texttt{Mixing Method++} is still an open problem, which could potentially be handled with the help of Lyapunov Analysis~\cite{wilson2016lyapunov,jin2018accelerated}.



\begin{acks}
AK acknowledges funding by the NSF (CCF-1907936). This work was partially done as MTT’s and CW’s class project for ``COMP545: Advanced Topics in Optimization,'' Rice University, Spring 2021. AK thanks Danny Carey for his percussion performance at the song ``Pneuma.''
\end{acks}

\bibliographystyle{ACM-Reference-Format}
\bibliography{biblio}


\begin{thebibliography}{64}


\ifx \showCODEN    \undefined \def \showCODEN     #1{\unskip}     \fi
\ifx \showDOI      \undefined \def \showDOI       #1{#1}\fi
\ifx \showISBNx    \undefined \def \showISBNx     #1{\unskip}     \fi
\ifx \showISBNxiii \undefined \def \showISBNxiii  #1{\unskip}     \fi
\ifx \showISSN     \undefined \def \showISSN      #1{\unskip}     \fi
\ifx \showLCCN     \undefined \def \showLCCN      #1{\unskip}     \fi
\ifx \shownote     \undefined \def \shownote      #1{#1}          \fi
\ifx \showarticletitle \undefined \def \showarticletitle #1{#1}   \fi
\ifx \showURL      \undefined \def \showURL       {\relax}        \fi
\providecommand\bibfield[2]{#2}
\providecommand\bibinfo[2]{#2}
\providecommand\natexlab[1]{#1}
\providecommand\showeprint[2][]{arXiv:#2}

\bibitem[\protect\citeauthoryear{Abbe}{Abbe}{2018}]%
        {abbe2018community}
\bibfield{author}{\bibinfo{person}{E. Abbe}.} \bibinfo{year}{2018}\natexlab{}.
\newblock \showarticletitle{Community Detection and Stochastic Block Models}.
\newblock \bibinfo{journal}{\emph{Foundations and Trends{\textregistered} in
  Communications and Information Theory}} \bibinfo{volume}{14},
  \bibinfo{number}{1-2} (\bibinfo{year}{2018}), \bibinfo{pages}{1--162}.
\newblock


\bibitem[\protect\citeauthoryear{Alizadeh, Haeberly, and Overton}{Alizadeh
  et~al\mbox{.}}{1997}]%
        {alizadeh1997complementarity}
\bibfield{author}{\bibinfo{person}{Farid Alizadeh},
  \bibinfo{person}{Jean-Pierre~A Haeberly}, {and} \bibinfo{person}{Michael~L
  Overton}.} \bibinfo{year}{1997}\natexlab{}.
\newblock \showarticletitle{Complementarity and nondegeneracy in semidefinite
  programming}.
\newblock \bibinfo{journal}{\emph{Mathematical programming}}
  \bibinfo{volume}{77}, \bibinfo{number}{1} (\bibinfo{year}{1997}),
  \bibinfo{pages}{111--128}.
\newblock


\bibitem[\protect\citeauthoryear{ApS}{ApS}{2019}]%
        {aps2019mosek}
\bibfield{author}{\bibinfo{person}{Mosek ApS}.}
  \bibinfo{year}{2019}\natexlab{}.
\newblock \showarticletitle{Mosek optimization toolbox for {MATLAB}}.
\newblock \bibinfo{journal}{\emph{User’s Guide and Reference Manual,
  version}}  \bibinfo{volume}{4} (\bibinfo{year}{2019}).
\newblock


\bibitem[\protect\citeauthoryear{Assran and Rabbat}{Assran and Rabbat}{2020}]%
        {assran2020convergence}
\bibfield{author}{\bibinfo{person}{M. Assran} {and} \bibinfo{person}{M.
  Rabbat}.} \bibinfo{year}{2020}\natexlab{}.
\newblock \showarticletitle{On the Convergence of {N}esterov's Accelerated
  Gradient Method in Stochastic Settings}.
\newblock \bibinfo{journal}{\emph{arXiv preprint arXiv:2002.12414}}
  (\bibinfo{year}{2020}).
\newblock


\bibitem[\protect\citeauthoryear{Barahona, Gr{\"o}tschel, J{\"u}nger, and
  Reinelt}{Barahona et~al\mbox{.}}{1988}]%
        {barahona1988application}
\bibfield{author}{\bibinfo{person}{F. Barahona}, \bibinfo{person}{M.
  Gr{\"o}tschel}, \bibinfo{person}{M. J{\"u}nger}, {and} \bibinfo{person}{G.
  Reinelt}.} \bibinfo{year}{1988}\natexlab{}.
\newblock \showarticletitle{An application of combinatorial optimization to
  statistical physics and circuit layout design}.
\newblock \bibinfo{journal}{\emph{Operations Research}} \bibinfo{volume}{36},
  \bibinfo{number}{3} (\bibinfo{year}{1988}), \bibinfo{pages}{493--513}.
\newblock


\bibitem[\protect\citeauthoryear{Barvinok}{Barvinok}{1995}]%
        {barvinok1995problems}
\bibfield{author}{\bibinfo{person}{Alexander~I. Barvinok}.}
  \bibinfo{year}{1995}\natexlab{}.
\newblock \showarticletitle{Problems of distance geometry and convex properties
  of quadratic maps}.
\newblock \bibinfo{journal}{\emph{Discrete \& Computational Geometry}}
  \bibinfo{volume}{13}, \bibinfo{number}{2} (\bibinfo{year}{1995}),
  \bibinfo{pages}{189--202}.
\newblock


\bibitem[\protect\citeauthoryear{Berg, Demirovi{\'{c}}, and Stuckey}{Berg
  et~al\mbox{.}}{2019}]%
        {loandra}
\bibfield{author}{\bibinfo{person}{Jeremias Berg}, \bibinfo{person}{Emir
  Demirovi{\'{c}}}, {and} \bibinfo{person}{Peter~J. Stuckey}.}
  \bibinfo{year}{2019}\natexlab{}.
\newblock \showarticletitle{Core-Boosted Linear Search for Incomplete MaxSAT}.
  In \bibinfo{booktitle}{\emph{Integration of Constraint Programming,
  Artificial Intelligence, and Operations Research}},
  \bibfield{editor}{\bibinfo{person}{Louis-Martin Rousseau} {and}
  \bibinfo{person}{Kostas Stergiou}} (Eds.). \bibinfo{publisher}{Springer
  International Publishing}, \bibinfo{address}{Cham}, \bibinfo{pages}{39--56}.
\newblock
\showISBNx{978-3-030-19212-9}


\bibitem[\protect\citeauthoryear{Bhojanapalli, Boumal, Jain, and
  Netrapalli}{Bhojanapalli et~al\mbox{.}}{2018}]%
        {bhojanapalli2018smoothed}
\bibfield{author}{\bibinfo{person}{S. Bhojanapalli}, \bibinfo{person}{N.
  Boumal}, \bibinfo{person}{P. Jain}, {and} \bibinfo{person}{P. Netrapalli}.}
  \bibinfo{year}{2018}\natexlab{}.
\newblock \showarticletitle{Smoothed analysis for low-rank solutions to
  semidefinite programs in quadratic penalty form}.
\newblock \bibinfo{journal}{\emph{arXiv preprint arXiv:1803.00186}}
  (\bibinfo{year}{2018}).
\newblock


\bibitem[\protect\citeauthoryear{Bhojanapalli, Kyrillidis, and
  Sanghavi}{Bhojanapalli et~al\mbox{.}}{2016}]%
        {bhojanapalli2016dropping}
\bibfield{author}{\bibinfo{person}{Srinadh Bhojanapalli},
  \bibinfo{person}{Anastasios Kyrillidis}, {and} \bibinfo{person}{Sujay
  Sanghavi}.} \bibinfo{year}{2016}\natexlab{}.
\newblock \showarticletitle{Dropping convexity for faster semi-definite
  optimization}. In \bibinfo{booktitle}{\emph{Conference on Learning Theory}}.
  \bibinfo{pages}{530--582}.
\newblock


\bibitem[\protect\citeauthoryear{Boumal}{Boumal}{2016}]%
        {boumal2016nonconvex}
\bibfield{author}{\bibinfo{person}{N. Boumal}.}
  \bibinfo{year}{2016}\natexlab{}.
\newblock \showarticletitle{Nonconvex phase synchronization}.
\newblock \bibinfo{journal}{\emph{SIAM Journal on Optimization}}
  \bibinfo{volume}{26}, \bibinfo{number}{4} (\bibinfo{year}{2016}),
  \bibinfo{pages}{2355--2377}.
\newblock


\bibitem[\protect\citeauthoryear{Burer and Monteiro}{Burer and
  Monteiro}{2003}]%
        {burer2003nonlinear}
\bibfield{author}{\bibinfo{person}{S. Burer} {and} \bibinfo{person}{R.
  Monteiro}.} \bibinfo{year}{2003}\natexlab{}.
\newblock \showarticletitle{A nonlinear programming algorithm for solving
  semidefinite programs via low-rank factorization}.
\newblock \bibinfo{journal}{\emph{Mathematical Programming}}
  \bibinfo{volume}{95}, \bibinfo{number}{2} (\bibinfo{year}{2003}),
  \bibinfo{pages}{329--357}.
\newblock


\bibitem[\protect\citeauthoryear{Chen, Wolfe, Li, and Kyrillidis}{Chen
  et~al\mbox{.}}{2019}]%
        {chen2019demon}
\bibfield{author}{\bibinfo{person}{John Chen}, \bibinfo{person}{Cameron Wolfe},
  \bibinfo{person}{Zhao Li}, {and} \bibinfo{person}{Anastasios Kyrillidis}.}
  \bibinfo{year}{2019}\natexlab{}.
\newblock \showarticletitle{Demon: Momentum Decay for Improved Neural Network
  Training}.
\newblock \bibinfo{journal}{\emph{arXiv preprint arXiv:1910.04952}}
  (\bibinfo{year}{2019}).
\newblock


\bibitem[\protect\citeauthoryear{Devolder, Glineur, and Nesterov}{Devolder
  et~al\mbox{.}}{2014}]%
        {devolder2014first}
\bibfield{author}{\bibinfo{person}{O. Devolder}, \bibinfo{person}{F. Glineur},
  {and} \bibinfo{person}{Y. Nesterov}.} \bibinfo{year}{2014}\natexlab{}.
\newblock \showarticletitle{First-order methods of smooth convex optimization
  with inexact oracle}.
\newblock \bibinfo{journal}{\emph{Mathematical Programming}}
  \bibinfo{volume}{146}, \bibinfo{number}{1-2} (\bibinfo{year}{2014}),
  \bibinfo{pages}{37--75}.
\newblock


\bibitem[\protect\citeauthoryear{Deza and Laurent}{Deza and Laurent}{1994}]%
        {deza1994applications}
\bibfield{author}{\bibinfo{person}{M. Deza} {and} \bibinfo{person}{M.
  Laurent}.} \bibinfo{year}{1994}\natexlab{}.
\newblock \showarticletitle{Applications of cut polyhedra - {II}}.
\newblock \bibinfo{journal}{\emph{J. Comput. Appl. Math.}}
  \bibinfo{volume}{55}, \bibinfo{number}{2} (\bibinfo{year}{1994}),
  \bibinfo{pages}{217--247}.
\newblock


\bibitem[\protect\citeauthoryear{{Erdogdu}, {Ozdaglar}, {Parrilo}, and
  {Denizcan Vanli}}{{Erdogdu} et~al\mbox{.}}{2018}]%
        {erdogdu2018convergence}
\bibfield{author}{\bibinfo{person}{Murat~A. {Erdogdu}}, \bibinfo{person}{Asuman
  {Ozdaglar}}, \bibinfo{person}{Pablo~A. {Parrilo}}, {and}
  \bibinfo{person}{Nuri {Denizcan Vanli}}.} \bibinfo{year}{2018}\natexlab{}.
\newblock \showarticletitle{{Convergence Rate of Block-Coordinate Maximization
  Burer-Monteiro Method for Solving Large SDPs}}.
\newblock \bibinfo{journal}{\emph{arXiv e-prints}}, Article
  \bibinfo{articleno}{arXiv:1807.04428} (\bibinfo{date}{July}
  \bibinfo{year}{2018}), \bibinfo{numpages}{arXiv:1807.04428}~pages.
\newblock
\showeprint[arxiv]{1807.04428}~[math.OC]


\bibitem[\protect\citeauthoryear{Frieze and Jerrum}{Frieze and Jerrum}{1995}]%
        {frieze1995improved}
\bibfield{author}{\bibinfo{person}{Alan Frieze} {and} \bibinfo{person}{Mark
  Jerrum}.} \bibinfo{year}{1995}\natexlab{}.
\newblock \showarticletitle{Improved approximation algorithms for {MAX}
  $k$-{CUT} and {MAX BISECTION}}. In \bibinfo{booktitle}{\emph{International
  Conference on Integer Programming and Combinatorial Optimization}}. Springer,
  \bibinfo{pages}{1--13}.
\newblock


\bibitem[\protect\citeauthoryear{G{\"a}rtner and Matousek}{G{\"a}rtner and
  Matousek}{2012}]%
        {gartner2012approximation}
\bibfield{author}{\bibinfo{person}{B. G{\"a}rtner} {and} \bibinfo{person}{J.
  Matousek}.} \bibinfo{year}{2012}\natexlab{}.
\newblock \bibinfo{booktitle}{\emph{Approximation algorithms and semidefinite
  programming}}.
\newblock \bibinfo{publisher}{Springer Science \& Business Media}.
\newblock


\bibitem[\protect\citeauthoryear{Ghadimi, Feyzmahdavian, and Johansson}{Ghadimi
  et~al\mbox{.}}{2015}]%
        {ghadimi2015global}
\bibfield{author}{\bibinfo{person}{E. Ghadimi}, \bibinfo{person}{H.
  Feyzmahdavian}, {and} \bibinfo{person}{M. Johansson}.}
  \bibinfo{year}{2015}\natexlab{}.
\newblock \showarticletitle{Global convergence of the heavy-ball method for
  convex optimization}. In \bibinfo{booktitle}{\emph{2015 European control
  conference (ECC)}}. IEEE, \bibinfo{pages}{310--315}.
\newblock


\bibitem[\protect\citeauthoryear{Gieseke, Pahikkala, and Igel}{Gieseke
  et~al\mbox{.}}{2013}]%
        {gieseke2013polynomial}
\bibfield{author}{\bibinfo{person}{F. Gieseke}, \bibinfo{person}{T. Pahikkala},
  {and} \bibinfo{person}{C. Igel}.} \bibinfo{year}{2013}\natexlab{}.
\newblock \showarticletitle{Polynomial runtime bounds for fixed-rank
  unsupervised least-squares classification}. In
  \bibinfo{booktitle}{\emph{Asian Conference on Machine Learning}}.
  \bibinfo{pages}{62--71}.
\newblock


\bibitem[\protect\citeauthoryear{Goemans and Williamson}{Goemans and
  Williamson}{1995}]%
        {goemans1995improved}
\bibfield{author}{\bibinfo{person}{M. Goemans} {and} \bibinfo{person}{D.
  Williamson}.} \bibinfo{year}{1995}\natexlab{}.
\newblock \showarticletitle{Improved approximation algorithms for maximum cut
  and satisfiability problems using semidefinite programming}.
\newblock \bibinfo{journal}{\emph{Journal of the ACM (JACM)}}
  \bibinfo{volume}{42}, \bibinfo{number}{6} (\bibinfo{year}{1995}),
  \bibinfo{pages}{1115--1145}.
\newblock


\bibitem[\protect\citeauthoryear{Goemans and Williamson}{Goemans and
  Williamson}{2004}]%
        {goemans2004approximation}
\bibfield{author}{\bibinfo{person}{Michel~X Goemans} {and}
  \bibinfo{person}{David~P Williamson}.} \bibinfo{year}{2004}\natexlab{}.
\newblock \showarticletitle{Approximation algorithms for MAX-3-CUT and other
  problems via complex semidefinite programming}.
\newblock \bibinfo{journal}{\emph{J. Comput. System Sci.}}
  \bibinfo{volume}{68}, \bibinfo{number}{2} (\bibinfo{year}{2004}),
  \bibinfo{pages}{442--470}.
\newblock


\bibitem[\protect\citeauthoryear{Goto, Tatsumura, and Dixon}{Goto
  et~al\mbox{.}}{2019}]%
        {goto2019combinatorial}
\bibfield{author}{\bibinfo{person}{Hayato Goto}, \bibinfo{person}{Kosuke
  Tatsumura}, {and} \bibinfo{person}{Alexander~R Dixon}.}
  \bibinfo{year}{2019}\natexlab{}.
\newblock \showarticletitle{Combinatorial optimization by simulating adiabatic
  bifurcat`ions in nonlinear {H}amiltonian systems}.
\newblock \bibinfo{journal}{\emph{Science advances}} \bibinfo{volume}{5},
  \bibinfo{number}{4} (\bibinfo{year}{2019}), \bibinfo{pages}{eaav2372}.
\newblock


\bibitem[\protect\citeauthoryear{Gr{\"o}tschel, Lov{\'a}sz, and
  Schrijver}{Gr{\"o}tschel et~al\mbox{.}}{2012}]%
        {grotschel2012geometric}
\bibfield{author}{\bibinfo{person}{M. Gr{\"o}tschel}, \bibinfo{person}{L.
  Lov{\'a}sz}, {and} \bibinfo{person}{A. Schrijver}.}
  \bibinfo{year}{2012}\natexlab{}.
\newblock \bibinfo{booktitle}{\emph{Geometric algorithms and combinatorial
  optimization}}. Vol.~\bibinfo{volume}{2}.
\newblock \bibinfo{publisher}{Springer Science \& Business Media}.
\newblock


\bibitem[\protect\citeauthoryear{Gurobi}{Gurobi}{2014}]%
        {optimization2014inc}
\bibfield{author}{\bibinfo{person}{Gurobi}.} \bibinfo{year}{2014}\natexlab{}.
\newblock \bibinfo{title}{{Inc.,“Gurobi} optimizer reference manual,”
  2015}.
\newblock
\newblock


\bibitem[\protect\citeauthoryear{Hajek, Wu, and Xu}{Hajek
  et~al\mbox{.}}{2016}]%
        {hajek2016achieving}
\bibfield{author}{\bibinfo{person}{B. Hajek}, \bibinfo{person}{Y. Wu}, {and}
  \bibinfo{person}{J. Xu}.} \bibinfo{year}{2016}\natexlab{}.
\newblock \showarticletitle{Achieving exact cluster recovery threshold via
  semidefinite programming}.
\newblock \bibinfo{journal}{\emph{IEEE Transactions on Information Theory}}
  \bibinfo{volume}{62}, \bibinfo{number}{5} (\bibinfo{year}{2016}),
  \bibinfo{pages}{2788--2797}.
\newblock


\bibitem[\protect\citeauthoryear{Hartmann}{Hartmann}{1996}]%
        {hartmann1996cluster}
\bibfield{author}{\bibinfo{person}{A. Hartmann}.}
  \bibinfo{year}{1996}\natexlab{}.
\newblock \showarticletitle{Cluster-exact approximation of spin glass
  groundstates}.
\newblock \bibinfo{journal}{\emph{Physica A: Statistical Mechanics and its
  Applications}} \bibinfo{volume}{224}, \bibinfo{number}{3-4}
  (\bibinfo{year}{1996}), \bibinfo{pages}{480--488}.
\newblock


\bibitem[\protect\citeauthoryear{Heath and Paulraj}{Heath and Paulraj}{1998}]%
        {heath1998simple}
\bibfield{author}{\bibinfo{person}{Robert~W Heath} {and}
  \bibinfo{person}{Arogyaswami Paulraj}.} \bibinfo{year}{1998}\natexlab{}.
\newblock \showarticletitle{A simple scheme for transmit diversity using
  partial channel feedback}. In \bibinfo{booktitle}{\emph{Conference Record of
  Thirty-Second Asilomar Conference on Signals, Systems and Computers (Cat. No.
  98CH36284)}}, Vol.~\bibinfo{volume}{2}. IEEE, \bibinfo{pages}{1073--1078}.
\newblock


\bibitem[\protect\citeauthoryear{Jin, Netrapalli, and Jordan}{Jin
  et~al\mbox{.}}{2018}]%
        {jin2018accelerated}
\bibfield{author}{\bibinfo{person}{Chi Jin}, \bibinfo{person}{Praneeth
  Netrapalli}, {and} \bibinfo{person}{Michael~I Jordan}.}
  \bibinfo{year}{2018}\natexlab{}.
\newblock \showarticletitle{Accelerated gradient descent escapes saddle points
  faster than gradient descent}. In \bibinfo{booktitle}{\emph{Conference On
  Learning Theory}}. \bibinfo{pages}{1042--1085}.
\newblock


\bibitem[\protect\citeauthoryear{Karp}{Karp}{1972}]%
        {karp1972reducibility}
\bibfield{author}{\bibinfo{person}{R. Karp}.} \bibinfo{year}{1972}\natexlab{}.
\newblock \showarticletitle{Reducibility among combinatorial problems}.
\newblock In \bibinfo{booktitle}{\emph{Complexity of computer computations}}.
  \bibinfo{publisher}{Springer}, \bibinfo{pages}{85--103}.
\newblock


\bibitem[\protect\citeauthoryear{Krislock, Malick, and Roupin}{Krislock
  et~al\mbox{.}}{2017}]%
        {krislock2017biqcrunch}
\bibfield{author}{\bibinfo{person}{N. Krislock}, \bibinfo{person}{J. Malick},
  {and} \bibinfo{person}{F. Roupin}.} \bibinfo{year}{2017}\natexlab{}.
\newblock \showarticletitle{Biq{C}runch: a semidefinite branch-and-bound method
  for solving binary quadratic problems}.
\newblock \bibinfo{journal}{\emph{ACM Transactions on Mathematical Software
  (TOMS)}} \bibinfo{volume}{43}, \bibinfo{number}{4} (\bibinfo{year}{2017}),
  \bibinfo{pages}{32}.
\newblock


\bibitem[\protect\citeauthoryear{Kyrillidis, Kalev, Park, Bhojanapalli,
  Caramanis, and Sanghavi}{Kyrillidis et~al\mbox{.}}{2018}]%
        {kyrillidis2018provable}
\bibfield{author}{\bibinfo{person}{Anastasios Kyrillidis},
  \bibinfo{person}{Amir Kalev}, \bibinfo{person}{Dohyung Park},
  \bibinfo{person}{Srinadh Bhojanapalli}, \bibinfo{person}{Constantine
  Caramanis}, {and} \bibinfo{person}{Sujay Sanghavi}.}
  \bibinfo{year}{2018}\natexlab{}.
\newblock \showarticletitle{Provable compressed sensing quantum state
  tomography via non-convex methods}.
\newblock \bibinfo{journal}{\emph{npj Quantum Information}}
  \bibinfo{volume}{4}, \bibinfo{number}{1} (\bibinfo{year}{2018}),
  \bibinfo{pages}{1--7}.
\newblock


\bibitem[\protect\citeauthoryear{Kyrillidis and Karystinos}{Kyrillidis and
  Karystinos}{2014}]%
        {kyrillidis2014fixed}
\bibfield{author}{\bibinfo{person}{Anastasios Kyrillidis} {and}
  \bibinfo{person}{George~N Karystinos}.} \bibinfo{year}{2014}\natexlab{}.
\newblock \showarticletitle{Fixed-rank Rayleigh quotient maximization by an
  MPSK sequence}.
\newblock \bibinfo{journal}{\emph{IEEE transactions on communications}}
  \bibinfo{volume}{62}, \bibinfo{number}{3} (\bibinfo{year}{2014}),
  \bibinfo{pages}{961--975}.
\newblock


\bibitem[\protect\citeauthoryear{Kyrillidis and Karystinos}{Kyrillidis and
  Karystinos}{2011}]%
        {kyrillidis2011rank}
\bibfield{author}{\bibinfo{person}{Anastasios~T Kyrillidis} {and}
  \bibinfo{person}{George~N Karystinos}.} \bibinfo{year}{2011}\natexlab{}.
\newblock \showarticletitle{Rank-deficient quadratic-form maximization over
  M-phase alphabet: Polynomial-complexity solvability and algorithmic
  developments}. In \bibinfo{booktitle}{\emph{2011 IEEE International
  Conference on Acoustics, Speech and Signal Processing (ICASSP)}}. IEEE,
  \bibinfo{pages}{3856--3859}.
\newblock


\bibitem[\protect\citeauthoryear{Lessard, Recht, and Packard}{Lessard
  et~al\mbox{.}}{2016}]%
        {lessard2016analysis}
\bibfield{author}{\bibinfo{person}{L. Lessard}, \bibinfo{person}{B. Recht},
  {and} \bibinfo{person}{A. Packard}.} \bibinfo{year}{2016}\natexlab{}.
\newblock \showarticletitle{Analysis and design of optimization algorithms via
  integral quadratic constraints}.
\newblock \bibinfo{journal}{\emph{SIAM Journal on Optimization}}
  \bibinfo{volume}{26}, \bibinfo{number}{1} (\bibinfo{year}{2016}),
  \bibinfo{pages}{57--95}.
\newblock


\bibitem[\protect\citeauthoryear{Liu, Yue, So, and Ma}{Liu
  et~al\mbox{.}}{2017}]%
        {liu2017discrete}
\bibfield{author}{\bibinfo{person}{Huikang Liu}, \bibinfo{person}{Man-Chung
  Yue}, \bibinfo{person}{Anthony Man-Cho So}, {and} \bibinfo{person}{Wing-Kin
  Ma}.} \bibinfo{year}{2017}\natexlab{}.
\newblock \showarticletitle{A discrete first-order method for large-scale MIMO
  detection with provable guarantees}. In \bibinfo{booktitle}{\emph{2017 IEEE
  18th International Workshop on Signal Processing Advances in Wireless
  Communications (SPAWC)}}. IEEE, \bibinfo{pages}{1--5}.
\newblock


\bibitem[\protect\citeauthoryear{Loizou and Richt{\'a}rik}{Loizou and
  Richt{\'a}rik}{2017}]%
        {loizou2017momentum}
\bibfield{author}{\bibinfo{person}{N. Loizou} {and} \bibinfo{person}{P.
  Richt{\'a}rik}.} \bibinfo{year}{2017}\natexlab{}.
\newblock \showarticletitle{Momentum and stochastic momentum for stochastic
  gradient, {N}ewton, proximal point and subspace descent methods}.
\newblock \bibinfo{journal}{\emph{arXiv preprint arXiv:1712.09677}}
  (\bibinfo{year}{2017}).
\newblock


\bibitem[\protect\citeauthoryear{Love, Heath, and Strohmer}{Love
  et~al\mbox{.}}{2003}]%
        {love2003grassmannian}
\bibfield{author}{\bibinfo{person}{David~J Love}, \bibinfo{person}{Robert~W
  Heath}, {and} \bibinfo{person}{Thomas Strohmer}.}
  \bibinfo{year}{2003}\natexlab{}.
\newblock \showarticletitle{Grassmannian beamforming for multiple-input
  multiple-output wireless systems}.
\newblock \bibinfo{journal}{\emph{IEEE transactions on information theory}}
  \bibinfo{volume}{49}, \bibinfo{number}{10} (\bibinfo{year}{2003}),
  \bibinfo{pages}{2735--2747}.
\newblock


\bibitem[\protect\citeauthoryear{Mart{\'\i}, Duarte, and Laguna}{Mart{\'\i}
  et~al\mbox{.}}{2009}]%
        {marti2009advanced}
\bibfield{author}{\bibinfo{person}{R. Mart{\'\i}}, \bibinfo{person}{A. Duarte},
  {and} \bibinfo{person}{M. Laguna}.} \bibinfo{year}{2009}\natexlab{}.
\newblock \showarticletitle{Advanced scatter search for the max-cut problem}.
\newblock \bibinfo{journal}{\emph{INFORMS Journal on Computing}}
  \bibinfo{volume}{21}, \bibinfo{number}{1} (\bibinfo{year}{2009}),
  \bibinfo{pages}{26--38}.
\newblock


\bibitem[\protect\citeauthoryear{Mei, Misiakiewicz, Montanari, and
  Oliveira}{Mei et~al\mbox{.}}{2017}]%
        {mei2017solving}
\bibfield{author}{\bibinfo{person}{S. Mei}, \bibinfo{person}{T. Misiakiewicz},
  \bibinfo{person}{A. Montanari}, {and} \bibinfo{person}{R. Oliveira}.}
  \bibinfo{year}{2017}\natexlab{}.
\newblock \showarticletitle{Solving {SDPs} for synchronization and {MaxCut}
  problems via the {G}rothendieck inequality}. In
  \bibinfo{booktitle}{\emph{Conference on Learning Theory}}.
  \bibinfo{pages}{1476--1515}.
\newblock


\bibitem[\protect\citeauthoryear{Mises and Pollaczek-Geiringer}{Mises and
  Pollaczek-Geiringer}{1929}]%
        {mises1929praktische}
\bibfield{author}{\bibinfo{person}{RV Mises} {and} \bibinfo{person}{Hilda
  Pollaczek-Geiringer}.} \bibinfo{year}{1929}\natexlab{}.
\newblock \showarticletitle{Praktische Verfahren der Gleichungsaufl{\"o}sung.}
\newblock \bibinfo{journal}{\emph{ZAMM-Journal of Applied Mathematics and
  Mechanics/Zeitschrift f{\"u}r Angewandte Mathematik und Mechanik}}
  \bibinfo{volume}{9}, \bibinfo{number}{2} (\bibinfo{year}{1929}),
  \bibinfo{pages}{152--164}.
\newblock


\bibitem[\protect\citeauthoryear{Mosek}{Mosek}{2015}]%
        {mosek2015mosek}
\bibfield{author}{\bibinfo{person}{ApS Mosek}.}
  \bibinfo{year}{2015}\natexlab{}.
\newblock \bibinfo{title}{The {MOSEK} optimization toolbox for {P}ython
  manual}.
\newblock
\newblock


\bibitem[\protect\citeauthoryear{Motedayen-Aval, Krishnamoorthy, and
  Anastasopoulos}{Motedayen-Aval et~al\mbox{.}}{2006}]%
        {motedayen2006optimal}
\bibfield{author}{\bibinfo{person}{Idin Motedayen-Aval},
  \bibinfo{person}{Arvind Krishnamoorthy}, {and} \bibinfo{person}{Achilleas
  Anastasopoulos}.} \bibinfo{year}{2006}\natexlab{}.
\newblock \showarticletitle{Optimal joint detection/estimation in fading
  channels with polynomial complexity}.
\newblock \bibinfo{journal}{\emph{IEEE transactions on information theory}}
  \bibinfo{volume}{53}, \bibinfo{number}{1} (\bibinfo{year}{2006}),
  \bibinfo{pages}{209--223}.
\newblock


\bibitem[\protect\citeauthoryear{Nesterov}{Nesterov}{2013}]%
        {nesterov2013introductory}
\bibfield{author}{\bibinfo{person}{Yurii Nesterov}.}
  \bibinfo{year}{2013}\natexlab{}.
\newblock \bibinfo{booktitle}{\emph{Introductory lectures on convex
  optimization: A basic course}}. Vol.~\bibinfo{volume}{87}.
\newblock \bibinfo{publisher}{Springer Science \& Business Media}.
\newblock


\bibitem[\protect\citeauthoryear{Nesterov and Nemirovskii}{Nesterov and
  Nemirovskii}{1989}]%
        {nesterov1989self}
\bibfield{author}{\bibinfo{person}{Y. Nesterov} {and} \bibinfo{person}{A.
  Nemirovskii}.} \bibinfo{year}{1989}\natexlab{}.
\newblock \bibinfo{booktitle}{\emph{Self-concordant functions and
  polynomial-time methods in convex programming}}.
\newblock \bibinfo{publisher}{USSR Academy of Sciences, Central Economic \&
  Mathematic Institute}.
\newblock


\bibitem[\protect\citeauthoryear{Nesterov and Nemirovskii}{Nesterov and
  Nemirovskii}{1994}]%
        {nesterov1994interior}
\bibfield{author}{\bibinfo{person}{Y. Nesterov} {and} \bibinfo{person}{A.
  Nemirovskii}.} \bibinfo{year}{1994}\natexlab{}.
\newblock \bibinfo{booktitle}{\emph{Interior-point polynomial algorithms in
  convex programming}}.
\newblock \bibinfo{publisher}{SIAM}.
\newblock


\bibitem[\protect\citeauthoryear{O’donoghue and Candes}{O’donoghue and
  Candes}{2015}]%
        {o2015adaptive}
\bibfield{author}{\bibinfo{person}{Brendan O’donoghue} {and}
  \bibinfo{person}{Emmanuel Candes}.} \bibinfo{year}{2015}\natexlab{}.
\newblock \showarticletitle{Adaptive restart for accelerated gradient schemes}.
\newblock \bibinfo{journal}{\emph{Foundations of computational mathematics}}
  \bibinfo{volume}{15}, \bibinfo{number}{3} (\bibinfo{year}{2015}),
  \bibinfo{pages}{715--732}.
\newblock


\bibitem[\protect\citeauthoryear{Pataki}{Pataki}{1998}]%
        {pataki1998rank}
\bibfield{author}{\bibinfo{person}{G{\'a}bor Pataki}.}
  \bibinfo{year}{1998}\natexlab{}.
\newblock \showarticletitle{On the rank of extreme matrices in semidefinite
  programs and the multiplicity of optimal eigenvalues}.
\newblock \bibinfo{journal}{\emph{Mathematics of operations research}}
  \bibinfo{volume}{23}, \bibinfo{number}{2} (\bibinfo{year}{1998}),
  \bibinfo{pages}{339--358}.
\newblock


\bibitem[\protect\citeauthoryear{Polyak}{Polyak}{1987}]%
        {polyak1987introduction}
\bibfield{author}{\bibinfo{person}{Boris~T Polyak}.}
  \bibinfo{year}{1987}\natexlab{}.
\newblock \showarticletitle{Introduction to optimization. optimization
  software}.
\newblock \bibinfo{journal}{\emph{Inc., Publications Division, New York}}
  \bibinfo{volume}{1} (\bibinfo{year}{1987}).
\newblock


\bibitem[\protect\citeauthoryear{Shi and Malik}{Shi and Malik}{2000}]%
        {shi2000normalized}
\bibfield{author}{\bibinfo{person}{J. Shi} {and} \bibinfo{person}{J. Malik}.}
  \bibinfo{year}{2000}\natexlab{}.
\newblock \showarticletitle{Normalized cuts and image segmentation}.
\newblock \bibinfo{journal}{\emph{IEEE Transactions on pattern analysis and
  machine intelligence}} \bibinfo{volume}{22}, \bibinfo{number}{8}
  (\bibinfo{year}{2000}), \bibinfo{pages}{888--905}.
\newblock


\bibitem[\protect\citeauthoryear{Singer}{Singer}{2011}]%
        {singer2011angular}
\bibfield{author}{\bibinfo{person}{A. Singer}.}
  \bibinfo{year}{2011}\natexlab{}.
\newblock \showarticletitle{Angular synchronization by eigenvectors and
  semidefinite programming}.
\newblock \bibinfo{journal}{\emph{Applied and computational harmonic analysis}}
  \bibinfo{volume}{30}, \bibinfo{number}{1} (\bibinfo{year}{2011}),
  \bibinfo{pages}{20--36}.
\newblock


\bibitem[\protect\citeauthoryear{So, Zhang, and Ye}{So et~al\mbox{.}}{2007}]%
        {so2007approximating}
\bibfield{author}{\bibinfo{person}{Anthony Man-Cho So}, \bibinfo{person}{Jiawei
  Zhang}, {and} \bibinfo{person}{Yinyu Ye}.} \bibinfo{year}{2007}\natexlab{}.
\newblock \showarticletitle{On approximating complex quadratic optimization
  problems via semidefinite programming relaxations}.
\newblock \bibinfo{journal}{\emph{Mathematical Programming}}
  \bibinfo{volume}{110}, \bibinfo{number}{1} (\bibinfo{year}{2007}),
  \bibinfo{pages}{93--110}.
\newblock


\bibitem[\protect\citeauthoryear{Sturm}{Sturm}{1999}]%
        {sturm1999using}
\bibfield{author}{\bibinfo{person}{Jos~F Sturm}.}
  \bibinfo{year}{1999}\natexlab{}.
\newblock \showarticletitle{Using {SeDuMi} 1.02, a {MATLAB} toolbox for
  optimization over symmetric cones}.
\newblock \bibinfo{journal}{\emph{Optimization methods and software}}
  \bibinfo{volume}{11}, \bibinfo{number}{1-4} (\bibinfo{year}{1999}),
  \bibinfo{pages}{625--653}.
\newblock


\bibitem[\protect\citeauthoryear{Tran-Dinh, Kyrillidis, and Cevher}{Tran-Dinh
  et~al\mbox{.}}{2014}]%
        {tran2014inexact}
\bibfield{author}{\bibinfo{person}{Q. Tran-Dinh}, \bibinfo{person}{A.
  Kyrillidis}, {and} \bibinfo{person}{V. Cevher}.}
  \bibinfo{year}{2014}\natexlab{}.
\newblock \showarticletitle{An inexact proximal path-following algorithm for
  constrained convex minimization}.
\newblock \bibinfo{journal}{\emph{SIAM Journal on Optimization}}
  \bibinfo{volume}{24}, \bibinfo{number}{4} (\bibinfo{year}{2014}),
  \bibinfo{pages}{1718--1745}.
\newblock


\bibitem[\protect\citeauthoryear{Tran-Dinh, Kyrillidis, and Cevher}{Tran-Dinh
  et~al\mbox{.}}{2016}]%
        {tran2016single}
\bibfield{author}{\bibinfo{person}{Q. Tran-Dinh}, \bibinfo{person}{A.
  Kyrillidis}, {and} \bibinfo{person}{V. Cevher}.}
  \bibinfo{year}{2016}\natexlab{}.
\newblock \showarticletitle{A single-phase, proximal path-following framework}.
\newblock \bibinfo{journal}{\emph{arXiv preprint arXiv:1603.01681}}
  (\bibinfo{year}{2016}).
\newblock


\bibitem[\protect\citeauthoryear{Veldt, Wirth, and Gleich}{Veldt
  et~al\mbox{.}}{2017}]%
        {veldt2017correlation}
\bibfield{author}{\bibinfo{person}{N. Veldt}, \bibinfo{person}{A. Wirth}, {and}
  \bibinfo{person}{D. Gleich}.} \bibinfo{year}{2017}\natexlab{}.
\newblock \showarticletitle{Correlation Clustering with Low-Rank Matrices}. In
  \bibinfo{booktitle}{\emph{Proceedings of the 26th International Conference on
  World Wide Web}}. International World Wide Web Conferences Steering
  Committee, \bibinfo{pages}{1025--1034}.
\newblock


\bibitem[\protect\citeauthoryear{Waldspurger, d’Aspremont, and
  Mallat}{Waldspurger et~al\mbox{.}}{2015}]%
        {waldspurger2015phase}
\bibfield{author}{\bibinfo{person}{I. Waldspurger}, \bibinfo{person}{A.
  d’Aspremont}, {and} \bibinfo{person}{S. Mallat}.}
  \bibinfo{year}{2015}\natexlab{}.
\newblock \showarticletitle{Phase recovery, {MAXCUT} and complex semidefinite
  programming}.
\newblock \bibinfo{journal}{\emph{Mathematical Programming}}
  \bibinfo{volume}{149}, \bibinfo{number}{1-2} (\bibinfo{year}{2015}),
  \bibinfo{pages}{47--81}.
\newblock


\bibitem[\protect\citeauthoryear{Wang, Nguyen, Bertozzi, Baraniuk, and
  Osher}{Wang et~al\mbox{.}}{2020}]%
        {wang2020scheduled}
\bibfield{author}{\bibinfo{person}{Bao Wang}, \bibinfo{person}{Tan~M Nguyen},
  \bibinfo{person}{Andrea~L Bertozzi}, \bibinfo{person}{Richard~G Baraniuk},
  {and} \bibinfo{person}{Stanley~J Osher}.} \bibinfo{year}{2020}\natexlab{}.
\newblock \showarticletitle{Scheduled restart momentum for accelerated
  stochastic gradient descent}.
\newblock \bibinfo{journal}{\emph{arXiv preprint arXiv:2002.10583}}
  (\bibinfo{year}{2020}).
\newblock


\bibitem[\protect\citeauthoryear{Wang, Jebara, and Chang}{Wang
  et~al\mbox{.}}{2013}]%
        {wang2013semi}
\bibfield{author}{\bibinfo{person}{J. Wang}, \bibinfo{person}{T. Jebara}, {and}
  \bibinfo{person}{S.-F. Chang}.} \bibinfo{year}{2013}\natexlab{}.
\newblock \showarticletitle{Semi-supervised learning using greedy {MaxCut}}.
\newblock \bibinfo{journal}{\emph{Journal of Machine Learning Research}}
  \bibinfo{volume}{14}, \bibinfo{number}{Mar} (\bibinfo{year}{2013}),
  \bibinfo{pages}{771--800}.
\newblock


\bibitem[\protect\citeauthoryear{Wang, Chang, and Kolter}{Wang
  et~al\mbox{.}}{2017}]%
        {wang2017mixing}
\bibfield{author}{\bibinfo{person}{P.-W. Wang}, \bibinfo{person}{W.-C. Chang},
  {and} \bibinfo{person}{Z. Kolter}.} \bibinfo{year}{2017}\natexlab{}.
\newblock \showarticletitle{The {M}ixing method: coordinate descent for
  low-rank semidefinite programming}.
\newblock \bibinfo{journal}{\emph{arXiv preprint arXiv:1706.00476}}
  (\bibinfo{year}{2017}).
\newblock


\bibitem[\protect\citeauthoryear{Wilson, Recht, and Jordan}{Wilson
  et~al\mbox{.}}{2016}]%
        {wilson2016lyapunov}
\bibfield{author}{\bibinfo{person}{Ashia~C Wilson}, \bibinfo{person}{Benjamin
  Recht}, {and} \bibinfo{person}{Michael~I Jordan}.}
  \bibinfo{year}{2016}\natexlab{}.
\newblock \showarticletitle{A {L}yapunov analysis of momentum methods in
  optimization}.
\newblock \bibinfo{journal}{\emph{arXiv preprint arXiv:1611.02635}}
  (\bibinfo{year}{2016}).
\newblock


\bibitem[\protect\citeauthoryear{Yang, Sun, and Toh}{Yang
  et~al\mbox{.}}{2015}]%
        {yang2015sdpnal}
\bibfield{author}{\bibinfo{person}{Liuqin Yang}, \bibinfo{person}{Defeng Sun},
  {and} \bibinfo{person}{Kim-Chuan Toh}.} \bibinfo{year}{2015}\natexlab{}.
\newblock \showarticletitle{{SDPNAL}: a majorized semismooth {Newton-CG}
  augmented Lagrangian method for semidefinite programming with nonnegative
  constraints}.
\newblock \bibinfo{journal}{\emph{Mathematical Programming Computation}}
  \bibinfo{volume}{7}, \bibinfo{number}{3} (\bibinfo{year}{2015}),
  \bibinfo{pages}{331--366}.
\newblock


\bibitem[\protect\citeauthoryear{Yurtsever, Tropp, Fercoq, Udell, and
  Cevher}{Yurtsever et~al\mbox{.}}{2019}]%
        {yurtsever2019scalable}
\bibfield{author}{\bibinfo{person}{Alp Yurtsever}, \bibinfo{person}{Joel~A
  Tropp}, \bibinfo{person}{Olivier Fercoq}, \bibinfo{person}{Madeleine Udell},
  {and} \bibinfo{person}{Volkan Cevher}.} \bibinfo{year}{2019}\natexlab{}.
\newblock \showarticletitle{Scalable Semidefinite Programming}.
\newblock \bibinfo{journal}{\emph{arXiv preprint arXiv:1912.02949}}
  (\bibinfo{year}{2019}).
\newblock


\bibitem[\protect\citeauthoryear{{Yurtsever}, {Tropp}, {Fercoq}, {Udell}, and
  {Cevher}}{{Yurtsever} et~al\mbox{.}}{2019}]%
        {CGAL}
\bibfield{author}{\bibinfo{person}{Alp {Yurtsever}}, \bibinfo{person}{Joel~A.
  {Tropp}}, \bibinfo{person}{Olivier {Fercoq}}, \bibinfo{person}{Madeleine
  {Udell}}, {and} \bibinfo{person}{Volkan {Cevher}}.}
  \bibinfo{year}{2019}\natexlab{}.
\newblock \showarticletitle{{Scalable Semidefinite Programming}}.
\newblock \bibinfo{journal}{\emph{arXiv e-prints}}, Article
  \bibinfo{articleno}{arXiv:1912.02949} (\bibinfo{date}{Dec}
  \bibinfo{year}{2019}), \bibinfo{numpages}{arXiv:1912.02949}~pages.
\newblock
\showeprint[arxiv]{1912.02949}~[math.OC]


\bibitem[\protect\citeauthoryear{Zhong and Boumal}{Zhong and Boumal}{2018}]%
        {zhong2018near}
\bibfield{author}{\bibinfo{person}{Y. Zhong} {and} \bibinfo{person}{N.
  Boumal}.} \bibinfo{year}{2018}\natexlab{}.
\newblock \showarticletitle{Near-optimal bounds for phase synchronization}.
\newblock \bibinfo{journal}{\emph{SIAM Journal on Optimization}}
  \bibinfo{volume}{28}, \bibinfo{number}{2} (\bibinfo{year}{2018}),
  \bibinfo{pages}{989--1016}.
\newblock


\end{thebibliography}

\newpage

\end{document}